\documentclass[12pt]{article}
\usepackage{amssymb}
\usepackage{amsmath}

\usepackage{hyperref}
\hypersetup{
    linktocpage = true,
    pdftitle = {Central extensions of groups of sections},
    pdfauthor = {Karl-Hermann Neeb and Christoph Wockel},
    bookmarksopen = true,
    bookmarksopenlevel = 1,
	unicode = true
  }

\input{liemacs.sty} 

\newcommand{\cyc}{{\rm cyc.}}

\newcommand{\Fl}{\ensuremath{\operatorname{Fl}}}
\newcommand{\tor}{\ensuremath{\operatorname{tor}}}
\newcommand{\wh}[1]{\ensuremath{\widehat{#1}}}
\newcommand{\wt}[1]{\ensuremath{\widetilde{#1}}}
\newcommand{\ol}[1]{\ensuremath{\overline{#1}}}
\newcommand{\cV}{\ensuremath{\mathcal{V}}}
\newcommand{\cK}{\ensuremath{\mathcal{K}}}
\newcommand{\from}{\ensuremath{\nobreak:\nobreak}}
\renewcommand{\to}{\ensuremath{\nobreak\rightarrow\nobreak}}
\newcommand{\pu} {\mathop{{\mathfrak{pu}}}\nolimits}
\newcommand{\D}{{\mathbb D}}
\newcommand{\cI}{{\mathcal I}}
\newcommand{\bV}{\V} 
\newcommand{\pdt}[1]{\frac{\partial{#1}}{\partial t}}
\newcommand{\pds}[1]{\frac{\partial{#1}}{\partial s}}

\begin{document}

\title{Central extensions of groups of sections} 

\author{Karl-Hermann Neeb and Christoph Wockel}

\maketitle

\begin{abstract} If $K$ is a Lie group and 
$q \: P \to M$ is a principal $K$-bundle over the 
compact manifold $M$, then any invariant symmetric $V$-valued bilinear form 
on the Lie algebra $\fk$ of $K$ defines a Lie algebra extension of the gauge algebra by 
a space of bundle-valued $1$-forms modulo exact $1$-forms. In the present 
paper we analyze the integrability of this extension to a Lie group 
extension for non-connected, possibly infinite-dimensional Lie
groups~$K$. If $K$ has finitely many connected components, we give a 
complete characterization of the integrable extensions. Our results 
on gauge groups are obtained by specialization of more general 
results on extensions of Lie groups of smooth sections of Lie group 
bundles. In this more general context we provide sufficient conditions 
for integrability in terms of data related only to the group~$K$. 
\end{abstract}

\section*{Introduction} 

Affine Kac--Moody groups and their 
Lie algebras play an interesting role in various fields of mathematics and 
mathematical physics, such as string theory and conformal field 
theory (cf.\ \cite{PS86} for the analytic theory of loop groups, 
\cite{Ka90} for the algebraic theory of Kac--Moody Lie algebras and  
\cite{Sch97} for connections to conformal field theory). 
For further connections to mathematical physics we refer to the monograph 
\cite{Mi89} which discusses various occurrences of Lie algebras of smooth 
maps in physical theories (see also \cite{Mu88}, \cite{DDS95}). 

From a geometric perspective, {\it affine Kac--Moody Lie groups} 
can be obtained from 
gauge groups $\Gau(P)$ of principal bundles $P$ over the circle $\bS^1$ 
whose fiber group is a simple compact Lie group $K$ 
by constructing a central extension and forming a semidirect product 
with a circle group corresponding to rigid rotations of the circle. 
Here the untwisted case corresponds to trivial bundles, where
 $\Gau(P) \cong C^\infty(\bS^1,K)$ is a loop group, 
and the twisted case corresponds to bundles which can be trivialized 
by a $2$- or $3$-fold covering of~$\bS^1$. 

In the present paper we address central extensions of gauge 
groups $\Gau(P)$ of more general bundles over a compact smooth manifold 
$M$, where the structure group $K$ may be an infinite-dimensional 
locally exponential Lie group. In particular, Banach--Lie groups and 
groups of smooth maps on compact manifolds are permitted. 
Since the gauge group $\Gau(P)$ is isomorphic to the group 
of smooth sections of the associated group bundle, defined by 
the conjugation action of $K$ on itself, it is natural to 
address central extensions of gauge groups and their Lie algebras 
in the more general context of groups of sections of bundles of Lie groups, 
resp., Lie algebras. 

In the following, $\cK$ always denotes a locally trivial Lie group 
bundle whose typical fiber $K$ is a locally exponential Lie group 
with Lie algebra $\fk = \L(K)$. 
Since we work with infinite-dimensional 
Lie algebras, we have to face the difficulty that, in general, 
the group $\Aut(\fk)$ does not carry a natural Lie group 
structure.\begin{footnote}{For Banach--Lie algebras, the group 
$\Aut(\fk)$ carries a natural Banach--Lie group structure with Lie algebra 
$\der(\fk)$, but if $\fk$ is not Banach, this need not be the case 
(cf.\ \cite{Mai63}).}
\end{footnote} Therefore it is natural to consider only Lie algebra bundles 
which are associated to some principal $H$-bundle $P$ 
with respect to a smooth action 
$\rho_\fk \: H \to \Aut(\fk)$ of a Lie group $H$ on $\fk$, 
i.e., for which the map $(h,x)\mapsto \rho_{\fk}(h)(x)$ is smooth.

Let $\fK$ be such a Lie algebra bundle. Then the smooth compact open topology 
turns the space $\Gamma\fK$ of its smooth sections into a locally 
convex topological Lie algebra. To construct $2$-cocycles 
on this algebra, we start with a continuous invariant symmetric bilinear 
map 
$$ \kappa \: \fk \times \fk \to V $$
with values in a locally convex $H$-module $V$ on which 
the identity component $H_0$ acts trivially.
The corresponding vector bundle $\V$ associated to $P$ is flat, 
so that we have a natural exterior derivative $\dd$ on $\V$-valued 
differential forms. If $V$ is finite-dimensional or 
$H$ acts on $V$ as a finite group, 
then $\dd(\Gamma \V)$ is a closed subspace of 
$\Omega^1(M,\V)$, so that the quotient 
$\oline\Omega^1(M,\V) := \Omega^1(M,\V)/\dd(\Gamma\V)$ 
inherits a natural Hausdorff topology (see the introduction to Section~\ref{sec:1}). 

We are interested in the cocycles on the Lie algebra $\Gamma\fK$ with values in 
the space $\oline\Omega^1(M,\V)$, given by 
\begin{equation}
  \label{eq:cocyc-def}
\omega_\kappa^\nabla(f,g) := [\kappa(f,\dd^\nabla g)]. 
\end{equation}
Here $\dd^\nabla$ is the covariant exterior differential 
on $\Omega^\bullet(M,\fK)$ induced by a principal connection 
$\nabla$ on~$P$.  For the special case of gauge algebras of 
principal bundles with connected compact structure group $K$, 
cocycles of this form have also been discussed briefly 
in \cite{LMNS98}. Clearly, \eqref{eq:cocyc-def} generalizes 
the well-known cocycles for 
Lie algebras of smooth maps, obtained from invariant bilinear forms 
and leading to universal central extensions of $\Gamma\fK$ if 
$\fK$ is trivial and $\fk$ is semisimple 
(cf.\ \cite{NW08a}, \cite{KL82}). Since we are presently far from 
a complete understanding of the variety of all central 
extensions of $\Gamma\fK$ or corresponding groups, it seems natural 
to study this class of cocycles first.
For other classes of cocycles, which are easier to handle, 
and their integrability we refer to \cite[Sect.~4]{Ne09} 
and \cite{Vi08}.

It seems quite likely that if $\fk$ is finite-dimensional semisimple, 
$\fK = \Ad(P)$ is the gauge bundle of a principal 
$K$-bundle and $\kappa$ is universal, 
then  the 
central extension of $\gau(P) \cong \Gamma\fK$ 
by $\oline\Omega^1(M,\V)$ defined by $\omega_\kappa^\nabla$ is 
universal. The analogous result for multiloop algebras 
has recently been obtained by E.~Neher (\cite[Thm.~2.13]{Neh07}, 
cf.\ also \cite[p.147]{PPS07}), so that one may be optimistic, at least 
if $M$ is a torus. 

Actually it is this class of examples that motivates the more 
complicated setting, where the group $H$ acts non-trivially on $V$. 
Already for twisted loop groups of real simple Lie algebras $\fk$, 
one is lead to non-connected structure 
groups and the universal target space $V(\fk)$ is a non-trivial module 
for $\Aut(\fk)$, on which $\Aut(\fk)_0$ acts 
trivially. 

The main goal of the present paper is to understand the 
integrability of the Lie algebra extension $\hat{\Gamma\fK}$ of $\Gamma\fK$ 
defined by the cocycle  $\omega :=\omega_{\kappa }:= \omega_\kappa^\nabla$ 
to a Lie group extension of the identity component 
of the Lie group $\Gamma\cK$ 
(cf.\ Appendix \ref{app:A} for the Lie group structure on this group). 
According to the general machinery for integrating 
central Lie algebra extensions described 
in \cite[Thm.~7.9]{Ne02a}, $\hat{\Gamma\fK}$ integrates to a Lie group 
extension of the identity component $(\Gamma\cK)_0$ if and only if 
the image $\Pi_\omega$ of the period homomorphism 
$$ \per_\omega \: \pi_2(\Gamma\cK) \to \oline\Omega^1(M,\V) $$
obtained by integration of the left invariant $2$-form on $\Gamma\cK$ 
defined by $\omega$ is discrete 
and the adjoint action of $\Gamma \fK$ on the central extension
$\hat{\Gamma \fK}$ integrates to an action of the 
corresponding connected Lie group $(\Gamma \cK)_0$ 
(cf.\ Appendix~\ref{app:C} for more details on these two conditions). 
Therefore our main task consists in verifying these two conditions, 
resp., in finding verifiable necessary and sufficient conditions for these 
conditions to be satisfied. 

To obtain information on the period group $\Pi_\omega$, 
it is natural to compose 
the cocycle with pullback maps 
$\gamma^* \: \oline\Omega^1(M,\V) \to \oline\Omega^1(\bS^1,\gamma^*\V), $
defined by smooth loops $\gamma \: \bS^1\to M$. To make this strategy 
work, we need quite detailed information on the special case 
$M = \bS^1$, for which $\Gamma\cK$ is a twisted loop group defined by 
some automorphism $\phi \in \Aut(K)$: 
$$ \Gamma\cK \cong C^\infty(\R,K)_\phi 
:= \{ f \in C^\infty(\R,K) \: (\forall t \in \R)\ 
f(t + 1) = \phi^{-1}(f(t))\}. $$

The structure of the paper is as follows. 
In Section~\ref{sec:1} we introduce the cocycles 
$\omega_\kappa^\nabla$ and discuss the dependence of their 
cohomology class on the connection $\nabla$. 
In particular, we show that they lift to $\Omega^1(M,\V)$-valued 
cocycles if $\kappa$ is exact in the sense that the 
$3$-cocycle $C(\kappa)(x,y,z) := \kappa([x,y],z)$ on $\fk$ is a coboundary 
(cf.\ \cite{Ne09}, \cite{NW08a}).  
In the end of this section we introduce an interesting class of bundles 
which are quite different from  adjoint bundles and illustrate many 
of the phenomena and difficulties we encounter in this paper. 

In Section~\ref{sec:2} we analyze the situation for the special case 
$M = \bS^1$, where $\Gamma\cK$ is a twisted loop group. 
It is a key observation that in this case the 
period map $\per_{\omega}$ is closely
related to the period map of the closed biinvariant $3$-form on $K$, 
determined by the Lie algebra $3$-cocycle $C(\kappa)$. 
To establish this relation, we need 
the connecting maps in the long exact 
homotopy sequence of the fibration defined by 
the evaluation map \mbox{$\ev_{0}^{K} \: C^\infty(\R,K)_\phi \to K,$} \mbox{$f \mapsto f(0)$}. 
Luckily, these connecting maps are given explicitly in
terms of the twist $\varphi$ and thus can be determined in concrete
examples. Having established the relation between $\per_{\omega}$
and $\per_{C(\kappa)}$, we use detailed knowledge on
$\per_{C(\kappa)}$ to derive conditions for the
discreteness of the image $\Pi_\omega$ of $\per_{\omega}$. 
In particular, we describe examples in which $\Pi_\omega$ is not discrete. 
For the case where $K$ is finite-dimensional, our results provide complete 
information, based on a detailed analysis of $\per_{C(\kappa_u)}$ for the 
universal invariant form $\kappa_u$ in Appendix~B. 

In Section~\ref{sec:3}, 
we turn to the integrability problem for a general compact manifold $M$. 
Our strategy is to 
compose with pullback homomorphisms $\gamma^{*} \: \Gamma\cK \to 
\Gamma(\gamma^*\cK)$, where $\gamma \: \bS^1 \to M$ is a smooth loop, 
and to determine under which
conditions the period homomorphism of the corresponding twisted 
loop group $\Gamma(\gamma^*\cK)$ only depends on the homotopy 
class of $\gamma$. If this condition is not satisfied, then our examples 
show that the period groups cannot be controlled in a reasonable way.
Fortunately, the latter condition is 
equivalent to the following requirement on the curvature $R(\theta)$ of the
principal connection $1$-form $\theta$ corresponding to $\nabla$ 
and the action $\L(\rho_{\fk})$: For each derivation
$D\in\im(\L(\rho_{\fk})\circ R(\theta))$, the periods of the $2$-cocycle
$\eta_{D}(x,y):=\kappa (x,Dy)$ have to vanish. This condition 
is formulated completely in terms of $K$ and it is always satisfied 
if $K$ is finite-dimensional because $\pi_2(K)$ vanishes in this case. 
If the curvature requirement is fulfilled, then 
$\Pi_\omega$ is contained in $H^{1}_{\rm dR}(M,\V)$, 
so that we can use integration maps 
$H^{1}_{\rm dR}(M,\V)\to V$ to reduce the
discreteness problem for $\Pi_\omega$ to bundles
over $\bS^{1}$. 

The second part of Section~\ref{sec:3} treats the lifting
problem for the important special case of gauge bundles 
$\fK = \Ad(P)$ and $\Gamma\cK = \Gau(P)$. 
In this case we even show that the action of the 
full automorphism group $\Aut(P)$ on $\Gamma\fK = \gau(P)$ 
lifts to an action on the central extension $\hat\gau(P)$, defined 
by $\omega$. We also give an integrability
criterion for this action to central extensions of the 
identity component $\Gau(P)_0$. 
Summarizing, we obtain for gauge bundles the following theorem: 

\begin{theorem} \label{thm:int-type1-gauge} 
If $\pi_0(K)$ is finite and $\dim K < \infty$, 
then the following are equivalent: 
\begin{description}
\item[\rm(1)] $\omega_{\kappa}$ integrates for each principal 
$K$-bundle $P$ over a compact manifold $M$ to a Lie group extension of 
$\Gau(P)_0$. 
\item[\rm(2)] $\omega_{\kappa}$ integrates for the trivial 
$K$-bundle $P = \bS^1 \times K$ over $M = \bS^1$ 
to a Lie group extension of $C^\infty(\bS^1,K)_0$. 
\item[\rm(3)] The image of $\per_{\omega_{\kappa}} \: \pi_3(K) \to V$ is discrete. 
\end{description}
These conditions are satisfied if $\kappa$ is the universal invariant 
symmetric bilinear form with values in $V(\fk)$. 
\end{theorem}

In order to increase the readability of the paper, we present some
background material in appendices. This comprises the Lie group
structure on groups of sections of Lie group bundles, a discussion of 
the universal invariant form for finite-dimensional Lie algebras, the
main results on integrating Lie algebra extensions to Lie group
extensions and some curvature issues for principal bundles, needed in 
Section~\ref{sec:3}.

\subsection*{Notation and basic concepts} 

A {\it Lie group} $G$ \index{infinite dimensional Lie group} 
is a group equipped with a 
smooth manifold structure modeled on a locally convex space 
for which the group multiplication and the 
inversion are smooth maps
(cf.\ \cite{Mil84}, \cite{Ne06} and \cite{GN09}). 
We write $\1 \in G$ for the identity element and 
$\lambda_g(x) = gx$, resp., $\rho_g(x) = xg$ for the left, resp.,
right multiplication on $G$. Then each $x \in T_\1(G)$ corresponds to
a unique left invariant vector field $x_l$ with 
$x_l(g) := T_\1(\lambda_g)x, g \in G.$
The space of left invariant vector fields is closed under the Lie
bracket of vector fields, hence inherits a Lie algebra structure. In
this sense we obtain on $T_\1(G)$ a continuous Lie bracket which
is uniquely determined by $[x,y]_l = [x_l, y_l]$ for $x,y \in T_\1(G)$. 
We write $\L(G) = \g$ for the so obtained locally convex Lie algebra and 
note that for morphisms $\phi \: G \to H$ of Lie groups we obtain 
with $\L(\phi) := T_\1(\phi)$ a functor from the category of Lie groups 
to the category of locally convex Lie algebras. 
We write $q_G \: \tilde G_0 \to G_0$ for the universal covering map of the identity 
component $G_0$ of $G$ and identify 
the discrete central subgroup $\ker q_G$ of $\tilde G_0$ with $\pi_1(G) \cong \pi_1(G_0)$. 

For a smooth map $f \: M \to G$ we define 
the (left) logarithmic derivative in $\Omega^1(M,\g)$ by 
$\delta(f)v_{m} := f(m)^{-1} \cdot T_m(f)v_{m}$, where $\cdot$ refers to the 
two-sided action of $G$ on its tangent bundle $TG$. 

In the following, we always write $I = [0,1]$ for the unit interval in $\R$. 
A Lie group $G$ is called {\it regular} if for each 
$\xi \in C^\infty(I,\g)$, the initial value problem 
$$ \gamma(0) = \1, \quad \gamma'(t) = \gamma(t)\cdot\xi(t) = T_\1(\lambda_{\gamma(t)})\xi(t) $$
has a solution $\gamma_\xi \in C^\infty(I,G)$, and the
evolution map 
$$ \evol_G \: C^\infty(I,\g) \to G, \quad \xi \mapsto \gamma_\xi(1) $$
is smooth (cf.\ \cite{Mil84}). 
For a locally convex space $E$, the regularity of the Lie group $(E,+)$ is equivalent 
to the Mackey completeness of $E$, i.e., to the existence of integrals 
of smooth curves $\gamma \: I \to E$. We also recall that for each regular 
Lie group $G$, its Lie algebra $\g$ is Mackey complete and 
that all Banach--Lie groups are regular (\cite{GN09}). 

A smooth map $\exp_G \: \g \to G$ is said to be an {\it exponential 
function} if for each $x \in \g$, the curve $\gamma_x(t) 
:= \exp_G(tx)$ is a 
homomorphism $\R \to G$ with $\gamma_x'(0) = x$. 
Presently, all known 
Lie groups modelled on complete locally convex spaces 
possess an exponential function. 
For Banach--Lie groups, its existence follows from the 
theory of ordinary differential equations in Banach spaces. 
A Lie group $G$ is called {\it locally exponential}, if 
it has an exponential function mapping an open $0$-neighborhood in 
$\g$ diffeomorphically onto an open neighborhood of $\1$ in $G$. 
For more details, we refer to Milnor's lecture notes \cite{Mil84}, 
the survey \cite{Ne06}, and the forthcoming monograph \cite{GN09}. 

If $q \: E \to B$ is a smooth fiber bundle, then we write 
$\Gamma E$ for its space of smooth sections. 

If $\g$ is a topological Lie algebra and $V$ a topological $\g$-module, 
we write $(C^\bullet(\g,V),\dd_\g)$ for the corresponding Lie algebra 
complex of continuous $V$-valued cochains (\cite{ChE48}).

\section{Central extensions of section algebras of Lie algebra bundles} 
\label{sec:1}

We now turn to the details and introduce our notation.  
We write $P(M,H,q_P)$ for an principal $H$-bundle over the smooth 
manifold $M$ with structure group $H$ and bundle projection 
$q_P \: P \to M$. To any such bundle $P$
and to any smooth action $\rho_{\fk}\from H\to \Aut(\fk)$%
, we associate the 
Lie algebra bundle $\fK$, 
which is the set $(P \times \fk)/H$ of $H$-orbits in $P \times \fk$ for the
action $h.(p,x) = (p.h^{-1}, \rho_\fk(h)x)$. We write $[(p,x)] := H.(p,x)$
for the elements of $\fK$ and $q_\fK \: \fK \to M, [(p,x)] \mapsto 
q_P(p)$ for the bundle projection.  

It is no loss of generality to assume that the bundle $P$ is connected. 
Indeed, if $P_1 \subeq P$ is a connected component, then 
$q(P_1) = M$ and 
$H_1 := \{ h \in H \: P_1.h = P_1\}$ is an open subgroup, so that 
$P_1$ is a principal $H_1$-bundle over $M$. Further, the canonical map 
$P_1 \times \fk \to \fK,(p,x) \mapsto [(p,x)]$ is surjective and 
induces a diffeomorphism $(P_1 \times \fk)/H_1 \cong \fK$. 
In the following we shall always assume that $P$ is connected. 
This implies that the connecting map $\delta_1 \: \pi_1(M) \to \pi_0(H)$ 
of the long exact homotopy sequence of $P$ is surjective. 

Further, let $V$ be a Fr\'echet $H$-module on which
the identity component $H_0$ acts trivially and 
$\rho_V \: H \to \GL(V)$ be the corresponding representation, so that 
$H_0 \subeq \ker \rho_V$ and $\rho_V$ factors through a 
representation $\oline\rho_V  \: \pi_0(H) \to \GL(V)$ 
of the discrete group $\pi_0(H)$. 
Accordingly, the associated vector bundle $\V := 
(P \times V)/H$ is flat. It is also associated 
via $\oline\rho_V$ to the 
squeezed bundle $P_0 := P/H_0$, which is a principal 
$\pi_0(H)$-bundle over $M$. 
Due to the flatness of $\V$, we  
have a natural exterior derivative $\dd$ on the space  
$\Omega^\bullet(M,\V)\cong \Omega^\bullet(P_0,V)^{\pi_0(H)}$ 
of $\V$-valued differential forms and 
we define $\oline\Omega^1(M,\V) := \Omega^1(M,\V)/\dd(\Gamma\V)$ 
and write its elements as $[\alpha]$, $\alpha \in \Omega^1(M,\V)$.
If $V$ is finite-dimensional, 
then $\dd \Gamma \bV$ is a closed subspace of 
the Fr\'echet space $\Omega^1(M,\bV)$, so that the quotient inherits 
a natural Hausdorff locally convex topology. 
In fact, in Lemma~\ref{lem:groupco} below we construct a 
continuous map $\Omega^1(M,\V) \to Z^1(\pi_1(M),V)$ (group cocycles 
with respect to the representation $\rho_{M}:=\rho_V \circ \delta_1$) 
and show that $\dd \Gamma\V$ is the inverse image of the space 
$B^1(\pi_1(M),V)$ of coboundaries which is finite-dimensional if 
$V$ is so, hence closed in the Fr\'echet space $Z^1(\pi_1(M),V)$. 
Therefore $\dd\Gamma\V$ is closed. 

If $\rho_V(H) = \rho_M(\pi_1(M))$ is finite, then 
$\hat M := \tilde M/\ker \rho_M$ is a finite covering manifold of $M$  
and $\D := \pi_1(M)/\ker \rho_M$ acts on $\hat M$  by deck transformations. 
We then have 
$\Omega^{\bullet}(M,\V) \cong \Omega^\bullet(\hat M,V)^\D$ 
and the finiteness of $\D$ implies that 
$\dd\Gamma\V \cong B^1_{\rm dR}(\hat M,V)^\D$, so that $\dd \Gamma\V$ 
is a closed subspace. 
We therefore assume in the following 
that either $\rho_V(H)$ is finite or that $V$ is finite-dimensional 
to ensure that $\oline\Omega^1(M,\V)$ carries a natural Fr\'echet space 
structure (cf.\ Remark~\ref{rem:hausdorff}).

Now let 
$\kappa \: \fk \times \fk \to V$
be an $H$-invariant continuous symmetric bilinear map which is also 
$\fk$-in\-va\-riant in the sense that 
$$ \kappa([x,y],z) = \kappa(x,[y,z])\quad \mbox{for} \quad x,y,z \in \fk. $$
The $H$-invariance of $\kappa$ implies that it defines a 
$C^\infty(M,\R)$-bilinear map 
$$ \Gamma\fK \times \Gamma \fK \to \Gamma\V,\quad
(f,g) \mapsto \kappa(f,g), \quad \kappa(f,g)(p) := \kappa(f(p), g(p)). $$
which defines
 a  $\Gamma\V$-valued invariant symmetric bilinear form on the
Lie algebra $\Gamma\fK$. 
To associate a Lie algebra $2$-cocycle to this data, we choose 
a principal connection $\nabla$ on the principal bundle $P$ and also write 
$\nabla$ for the associated connections on the vector bundles 
$\fK$ and $\V$. Since $H$ acts by automorphisms on $\fk$, its Lie algebra 
$\fh$ acts by derivations, which implies that the connection 
$\nabla$ on $\fK$ is a {\it Lie connection}, i.e., 
\begin{equation}
  \label{eq:1.1}
 \nabla_X [f,g] = [\nabla_X f,g] + [f,\nabla_X g] 
\quad \mbox { for }\quad  \quad X \in {\cal V}(M), 
f,g \in \Gamma\fK 
\end{equation}
(cf.\ \cite{Ma05} for more details on Lie connections 
on Lie algebra bundles). The $H$-invariance of $\kappa$ and 
the fact that its Lie algebra $\fh$ acts trivially on $V$ imply that 
\begin{equation}
  \label{eq:1.2}
\dd\big(\kappa(f,g)(X)\big) = \kappa(\nabla_X f,g) + \kappa(f,\nabla_X g)\quad 
\mbox{ for } \quad X \in {\cal V}(M), 
f,g \in \Gamma\fK.   
\end{equation}
In the following we write $\dd^\nabla f$ for the $\fK$-valued $1$-form 
defined for $f \in \Gamma\fK$ by 
$(\dd^\nabla f)(X) := \nabla_X f$ for $X \in {\cal V}(M)$. 
In the realization of $\Gamma\fK$ as $C^\infty(P,\fk)^H$, we have 
for $X \in {\cal V}(P)$: 
$$ (\dd^\nabla f)(X) = \dd f(X) + \theta(X)f, $$
where $\theta \in \Omega^1(P,\fh)$ is the principal connection $1$-form 
corresponding to the connection~$\nabla$.

\begin{proposition} \label{prop:1.1}  The prescription 
$$ \omega(f,g) := \omega_\kappa^\nabla(f,g) := [\kappa(f,\dd^\nabla g)] $$
defines a Lie algebra cocycle on $\Gamma \fK$ with values in the trivial 
$\Gamma \fK$-module $\oline\Omega^1(M,\V)$. 
If $\nabla' = \nabla + \beta$, $\beta \in \Omega^1(M,\Ad(P))$, 
is another principal connection for which there exists some 
$\gamma \in \Omega^1(M,\fK)$ with 
\begin{equation}
  \label{eq:1.3} 
\L(\rho_\fk)\circ \beta(X) = \ad(\gamma(X)) \in \End(\Gamma\fK) 
\quad \mbox{ for } \quad X \in {\cal V}(M), 
\end{equation}
then the corresponding 
cocycle $\omega'$ differs from $\omega$ by a coboundary.  
\end{proposition}

\begin{proof} From \eqref{eq:1.2} we get 
$\dd\big(\kappa(f,g)\big) = \kappa(\dd^\nabla f,g) + \kappa(f,\dd^\nabla g)$, 
so that $\omega$ is alternating. 
In view of \eqref{eq:1.1} and \eqref{eq:1.2}, we further have 
\begin{eqnarray*}
\dd\big(\kappa([f,g],h)\big)
&=& \kappa(\dd^\nabla[f,g], h) + \kappa([f,g], \dd^\nabla h) \\ 
&=& \kappa([\dd^\nabla f,g], h) + \kappa([f,\dd^\nabla g],h) +  \kappa([f,g], 
\dd^\nabla h) \\ 
&=& \kappa([g, h], \dd^\nabla f) + \kappa([h,f],\dd^\nabla g) 
+  \kappa([f,g], \dd^\nabla h), 
\end{eqnarray*}
showing that $\omega$ is a $2$-cocycle. 

If $\nabla$ is replaced by $\nabla' = \nabla + \beta$ and 
\eqref{eq:1.3} is satisfied, 
then 
$$ \kappa(f,\dd^{\nabla'} g) 
=  \kappa(f,\dd^\nabla g) + \kappa(f,[\gamma,g])
=  \kappa(f,\dd^\nabla g) + \kappa(\gamma,[g,f]) $$
implies that 
$\omega' - \omega = \dd_{\Gamma\fK}([\kappa(\gamma,\cdot)]),$
where $\kappa(\gamma,\cdot)$ is an $\Omega^1(M,\V)$-valued 
linear map on $\Gamma\fK$. 
\end{proof}

\begin{remark} \label{rem:1.2} 
Since the space $\oline\Omega^1(M,\V)$ is a quotient 
of the space $\Omega^1(M,\bV)$ of $\V$-valued $1$-forms, it is 
natural to ask for the existence of $\Omega^1(M,\V)$-valued 
cocycles on $\Gamma\fK$ lifting $\omega_\kappa^\nabla$. To 
see when such cocycles exist, we consider the continuous bilinear 
map 
$$ \tilde\omega(f,g) := \kappa(f,\dd^\nabla g) - \kappa(g,\dd^\nabla f), $$
which is an alternating lift of $2\omega_\kappa^\nabla$. Its Lie algebra 
differential is 
\begin{align*}
(\dd_{\Gamma\fK}\tilde\omega)(f,g,h) 
&= -\sum_\cyc \big(\kappa([f,g],\dd^\nabla h) - \kappa(h,\dd^\nabla [f,g])\big)  = \sum_\cyc \kappa([f,g],\dd^\nabla h) \\
&= \dd(\kappa([f,g],h)), 
\end{align*}
as we see with similar calculations as in the proof 
Proposition~\ref{prop:1.1}. 

For the trivial $\fk$-module $V$, we write $\Sym^2(\fk,V)^\fk$ for the 
space of $V$-valued symmetric invariant bilinear forms, and recall the 
{\it Cartan map} 
$$ C \: \Sym^2(\fk,V)^\fk \to Z^3(\fk,V), \quad 
C(\kappa)(x,y,z) := \kappa([x,y],z). $$
We say that  $\kappa$ is {\it exact} if $C(\kappa)$ is a coboundary.
If $C(\kappa) = \dd_\fk \eta$ for some 
$\eta \in C^2(\fk,V)$, then 
$$\dd(\kappa([f,g],h)) = \dd((\dd_\fk\eta)(f,g,h))
= - \sum_\cyc \dd(\eta([f,g],h)), $$
so that 
$$ \omega_{\kappa,\eta}(f,g) :=  
\kappa(f,\dd^\nabla g) - \kappa(g,\dd^\nabla f) - \dd(\eta(f,g)) $$
is an $\Omega^1(M,\V)$-valued $2$-cocycle on $\Gamma\fK$ lifting 
$2\omega_\kappa^\nabla$ (cf.\ \cite[Sect.~2]{Ne09}). 
\end{remark}

\begin{remark} \label{rem:1.3} 
If $\beta \in \Omega^1(M,\Ad(P))$ is a bundle-valued 
$1$-form, then we obtain for each $X \in {\cal V}(M)$ 
a derivation $\beta_\fk(X) := \rho_\fk \circ \beta(X)$ of $\Gamma\fK$ 
and this derivation preserves the symmetric bilinear $\Gamma\V$-valued 
map $(f,g) \mapsto \kappa(f,g)$,  so that 
$$ \eta_\beta(f,g) := \kappa(f,\beta g) $$
defines an $\Omega^1(M,\V)$-valued $2$-cocycle on $\Gamma\fK$. 
For $\nabla' = \nabla + \beta$, we now have 
$$ \omega' - \omega = q_\Omega\circ \eta_\beta, $$
where $q_\Omega \: \Omega^1(M,\V) \to \oline\Omega^1(M,\V)$ denotes 
the quotient map. This argument shows that the dependence of the 
cohomology class $[\omega_\kappa^\nabla]$ on $\nabla$ is described 
by elements of
$H^2(\Gamma\fK, \Omega^1(M,\V))$. 

We may also consider $\eta_\beta$ as a bundle map $\fK \times \fK \to 
\Hom(TM,\V)$, which implies that $\eta_\beta$ can also be used to define a 
central extension of Lie algebroids (cf.\ \cite{Ma05}). 
\end{remark}

\begin{example} \label{ex:1.4} 
Of particular importance is the special case where 
$K := H$ is a Lie group with Lie algebra $\fk$ and 
$\rho_\fk \: K \to \Aut(\fk)$ is the adjoint action of $K$. 
Then $\fK = \Ad(P)$ is the adjoint bundle of the principal $K$-bundle $P$ 
over $M$, $\Gamma\fK \cong \gau(P)$, and we have the Lie algebra extension 
$$ \0 \to \gau(P) \cong \Gamma\fK\into \aut(P) = {\cal V}(P)^K 
\onto {\cal V}(M) \to \0. $$
Furthermore, the space $\Gamma\V$ of smooth sections of the flat 
vector bundle $\bV$ 
carries a natural ${\cal V}(M)$-module 
structure which we may pull back to an $\aut(P)$-module structure 
for which the ideal $\gau(P)$ acts trivially. 

Since $\L(\rho_\fk) = \ad$, Proposition~\ref{prop:1.1} implies that 
in this situation the cohomology class 
$$ [\omega] = [\omega_\kappa^\nabla] \in H^2(\gau(P), \oline\Omega^1(M,\V)) $$
does not depend on the choice of the principal connection in $P$. 
\end{example}

\begin{remark} If $\kappa_u \: \fk \times \fk \to V(\fk)$ is the
universal continuous invariant symmetric bilinear form on $\fk$ (cf.\
\cite{MN03} and Appendix~\ref{app:B} below) 
and $K$ is a Lie group with Lie algebra $\fk$, then the
universality of $\kappa_u$ implies that $K$ acts 
naturally on $V(\fk)$, and since
$\fk$ acts trivially on $V(\fk)$, the identity component $K_0$ acts
trivially (\cite{GN09}; \cite[Rem.~II.3.7]{Ne06}). 
This implies that the universal form 
$\kappa_u$ satisfies all assumptions required for our construction.
For a detailed analysis of $\kappa_u$ and the period map of the corresponding 
closed $3$-form on $K$, we refer to Appendix~\ref{app:B}.
\end{remark}

The aim of this paper is to determine under which circumstances the
Lie algebra extension defined by the cocycle $\omega$ from
Proposition~\ref{prop:1.1} integrates to an extension of Lie groups. The
natural setting for this question is the case, where the action
$\rho_{\fk}$ is induced by a smooth action $\rho_{K}\from
H\to\Aut(K)$, i.e., $K$ is a Lie group with Lie algebra $\fk$
and we have $\rho_{\fk}=\L(\rho_{K})$. 
If $K$ is locally exponential, then the group of sections
$\Gamma\cK$ of the adjoint Lie group bundle $\cK :=(P\times K)/H$ has
a natural Lie group structure with $\L(\Gamma\cK) \cong \Gamma \fK$ (cf.\
Appendix~\ref{app:A}). We therefore want to integrate our Lie algebra 
extension to the identity component $(\Gamma \cK)_0$ of this group. 

From \cite{Ne02a} (cf.\ Appendix~\ref{app:C}) 
we know that the Lie algebra cocycle 
$\omega$ defines a period map 
$$ \per_\omega \: \pi_2(\Gamma{\cal K}) \to \oline\Omega^1(M,\V), $$
and a necessary condition for the existence of a Lie group extension
integrating $\omega$ is that the image $\Pi_\omega$ of the period map,
the {\it period group}, is discrete (\cite{Ne02a}, Theorem~VII.9).
To obtain information on this period group, 
our strategy is first to take a closer look at the case $M = \bS^1$ and 
then to use this case to treat more general situations. 
The much simpler 
case of trivial bundles has been treated in a similar fashion in 
\cite{MN03}. 

\subsection*{A class of examples} 

\begin{example}\label{ex:Diff(N)-bundles}
Let $\pi \: Q\to M$
be a compact locally trivial smooth bundle with (compact) fiber
$N$. Then $Q$ is associated to the principal $H:=\Diff(N)$-bundle
$P$ with fiber $P_{m}:=\Diff(N,Q_{m})$ with the canonical $H$
action by composition. 
For any locally convex Lie group $G$, we have a canonical
$H$-action on $K := C^{\infty}(N,G)$ by $(\varphi ,\gamma)\mapsto 
\gamma \circ \varphi^{-1}$ whose smoothness follows from the smoothness 
of the action of $\Diff(N)$ on $N$ and the smoothness of the evaluation 
map of $K$ (cf.\ \cite{NW08b}, Lemma~A.2) and we
thus obtain an associated Lie group bundle
$\cK:=P\times_{H}C^{\infty}(N,G)$. The sections of the 
Lie group bundle $\Gamma \cK$ may be identified with the set 
$C^{\infty}(P,K)^{H}$ of $H$-equivariant smooth functions $P \to K$. 
\end{example}

\begin{proposition} If $G$ is locally exponential, then the map 
$$ s\from C^{\infty}(Q,G)\to \Gamma \cK, \quad s_{f}(p)=f\circ p,\quad 
p\in P_{m}=\Diff(N,P_{m})$$
is an isomorphism of Lie groups. 
\end{proposition}

\begin{proof} In local
coordinates one easily checks that $s$ actually is an 
isomorphism of abstract 
groups compatible with the smooth compact open topology, so that it 
actually is an isomorphism of topological groups.

If, in addition, $G$ is assumed to be locally exponential, then 
$\Gamma \cK$ and $C^\infty(Q,G)$ inherit this property 
(Theorem~\ref{thm:liestruc}), and now the general theory of locally 
exponential Lie groups (\cite[Thm.~IV.1.18]{Ne06}, \cite{GN09}) implies that 
the topological isomorphism between these groups actually 
is a diffeomorphism, hence an isomorphism of Lie groups. 
\end{proof}

\begin{remark} \label{rem:princbun-ex}
If the bundle $\pi \: Q \to M$ 
in Example~\ref{ex:Diff(N)-bundles} is an principal 
$H$-bundle for some compact group $H$, then the structure group 
can be reduced from the infinite-dimensional Lie group $\Diff(H)$ 
to the compact subgroup $H$, because the transition functions 
of the bundle charts have values in the group of left multiplications 
of~$H$. We then obtain an isomorphism of Lie groups
$$ s \: C^\infty(P,G) \to \Gamma\cK \cong C^\infty(P,K)^H, \quad 
s_f(p)(h) := f(p.h). $$
We thus associate to each principal $H$-bundle $Q$ a Lie group bundle 
$\cK$ with fiber $K = C^\infty(H,G)$. This construction is particularly 
interesting for $H = \bT$. Then $P$ is a circle bundle and 
$K = C^\infty(\T,G)$ is the loop group of $G$. 
\end{remark}

\section{Lie group bundles over the circle} 
\label{sec:2}

Throughout this section we consider the special case 
$M = \bS^1$ and assume that the Lie group $K$ is regular. 
Then every $K$-Lie group bundle over $\bS^1$ is flat, hence determined 
by its holonomy $\phi \in \Aut(K)$. Conversely, every automorphism 
$\phi \in \Aut(K)$ leads to a Lie group bundle 
${\cal K}_{\phi} = \R \times_{\phi} K$ over $\bS^1$ with holonomy $\phi$. 
Indeed, $\cK_\phi$ 
is the Lie group bundle associated to the universal 
covering $q_{\bS^1} \: \R \to \bS^1 \cong \R/\Z$ by the action of 
$H := \Z \cong \pi_1(\bS^1)$ on $K$ defined by $\phi$. 
The smooth sections of ${\cal K}_{\phi}$ 
correspond to \textit{twisted loops}: 
$$ \Gamma\cK_{\phi} \cong C^\infty(\R,K)_{\phi} := 
\{ f \in C^\infty(\R,K) \: (\forall t \in \R)\ 
f(t + 1) = \phi^{-1}(f(t))\}. $$
From now on we identify $\Gamma \cK_{\varphi}$ with 
$C^\infty(\R,K)_{\phi}$ and  write  
$$\ev_0^K \: C^\infty(\R,K)_{\phi} \to K, \quad f \mapsto f(0),$$ 
for the evaluation homomorphism in $0$. On the Lie algebra level, we
similarly get with $\phi_\fk = \L(\phi)$ a Lie algebra bundle
$\fK_{\varphi_{\fk}}$ with 
$$ \Gamma\fK_{\phi_\fk} \cong C^\infty(\R,\fk)_{\phi_\fk} 
:= \{ f \in C^\infty(\R,\fk) \: (\forall t \in \R)\ 
f(t + 1) = \phi_\fk^{-1}(f(t))\} $$
and $\L(\ev_0^K) = \ev_0^\fk$. 

\subsection{On the topology of twisted loop groups} \label{subsec:2.1}

In Section~\ref{sec:3} we shall reduce the calculation 
of the period groups for $\omega_\kappa^\nabla$ 
essentially to the case $M = \bS^1$, so that we need detailed information 
on the second homotopy group of twisted loop groups. 
A central tool is a simple description of the 
connecting maps in the long exact homotopy sequence defined 
by the evaluation map $\ev_0^K$ 
for a twisted loop group $C^\infty(\R,K)_\phi$, 
which is based on the fact that the passage from smooth to continuous 
twisted loops is a weak homotopy equivalence. 

\begin{lemma} The image of the evaluation homomorphism 
$\ev_0^K$ is the open subgroup 
$$ K^{[\phi]} := \{ k \in K \: \phi(k)k^{-1} \in K_0\} 
=  \{ k \in K \: kK_0 \in \pi_0(K)^\phi\}.$$
\end{lemma}

\begin{proof} For $f \in C^\infty(\R,K)_\phi$, we have 
$f(1) = \phi^{-1}(f(0))\in f(0)K_0$, so that 
the image of $\ev_0^K$ is contained in 
$K^{[\phi]}$. 
If, conversely, $k \in K^{[\phi]}$, then there exists a smooth 
curve $\alpha : [0,1] \to K$ with $\alpha(0) = k$, $\alpha(1) 
=\phi^{-1}(k)$ such that $\alpha$ is locally constant near $0$ and $1$. Then 
$f(n + t) := \phi^{-n}(\alpha(t))$ for $t \in [0,1], n \in \Z$, 
defines a section of ${\cal K}_\phi$ with $f(0) = k$. 
\end{proof}

\begin{lemma} \label{lem:loc-triv} The Lie group homomorphism $\ev_0^K$ 
has smooth local sections, hence defines a Lie group extension 
of $K^{[\phi]}$ by $C^\infty_*(\R,K)_\phi := \ker \ev_0^K$.  
\end{lemma}

\begin{proof} That $\ev_0^K$ has smooth local sections can be seen as follows. 
Let $(\psi,U)$ be a chart of $K$, centered in $\1$ for which 
$\psi(U)$ is convex. Let further $h \: [0,1] \to \R$ be a smooth 
function with $h(0) = 0$ and $h(1) = 1$ which is constant in 
$[0,\eps]$ and $[1-\eps,1]$. 
Let $W \subeq U \cap \phi(U)$ be a $\1$-neighborhood in $K$. 
For $k \in W$ we then consider the smooth curve 
$$ \gamma_k \: [0,1] \to K, \quad 
\gamma_k(t) := \psi^{-1}\big((1-h(t))\psi(k) + h(t)\psi(\phi^{-1}(k))\big). $$
Then $\gamma_k$ is constant near $0$ and $1$, 
$\gamma_k(0) = k$, and $\gamma_k(1) = \phi^{-1}(k)$. 
We extend $\gamma_k$ smoothly to $\R$ in such a way that it defines an 
element on $C^\infty(\R,K)_\phi$. Then the smoothness of the map 
$W \to C^\infty(\R,K)_\phi, k \mapsto \gamma_k$
follows from the smoothness of the corresponding map 
$W \times \R \to K, (k,t) \mapsto \gamma_k(t),$ 
which in turn follows from its smoothness on each subset $W \times ]n-\eps,n+1+\eps[$ (cf.\ \cite{GN09}).  
\end{proof} 

\begin{proposition} \label{prop:3.3} The inclusion 
$C^\infty(\R,K)_\phi \into C(\R,K)_\phi$ of the smooth 
twisted loop group into the continuous twisted loop group 
is a weak homotopy equivalence. 
\end{proposition}

\begin{proof} Let $H := K \rtimes_\phi \Z$, 
where the action of $\Z$ on $K$ is 
  defined by $n.k := \phi^n(k)$ for $n \in \Z$ and $k \in K$. 
We write $P_\phi$ for the principal $H$-bundle over $\bS^1$ 
with holonomy $(\1,1) \in H$. Then 
$$ \Gau(P_\phi) 
\cong \{ f \in C^\infty(\R,H) \: (\forall t \in \R)\ 
f(t + 1) = (\1,1)f(t)(\1,-1)\}, $$
and this group contains the twisted loop group $C^\infty(\R,K)_\phi$ 
as an open subgroup. According to Prop.~1.20 in \cite{Wo07a}, the 
inclusion of the smooth gauge group $\Gau(P_\phi)$ into the 
group $\Gau^c(P_\phi)$ of continuous gauge transformations 
is a weak homotopy equivalence, and this property is inherited 
by the open subgroups  of $K$-valued twisted loops. 
\end{proof}

\begin{corollary} \label{cor:3.4}The inclusion 
$C^\infty_*(\R,K)_\phi \into C_*(\R,K)_\phi$ of the smooth 
based twisted loop group into the continuous based twisted loop group 
is a weak homotopy equivalence. 
\end{corollary}

\begin{proof} In view of Lemma~\ref{lem:loc-triv}, the evaluation 
$\ev_0^K$ defines a smoothly locally trivial fiber bundle 
$C^\infty(\R,K)_\phi \to K^{[\phi]}$, and a similar (even simpler) 
argument shows that the same holds for the continuous twisted loop group. 
Since  $\id_{K^{[\phi]}}$ and the inclusion 
$C^\infty(\R,K)_\phi \into C(\R,K)_\phi$ are weak homotopy 
equivalences, the 5-Lemma implies that the same holds for the inclusion   
$C^\infty_*(\R,K)_\phi \into C_*(\R,K)_\phi$ of the fibers 
(cf.\ Prop.~A.8 in \cite{Ne02c}). 
\end{proof}

We have already determined the image of $\ev_0^K$, showing that the 
long exact homotopy sequence ends with 
$$ \ldots\to \pi_1(K) \ssmapright{\delta_1} \pi_0(C^\infty_*(\R,K)_\phi) 
\to \pi_0(C^\infty(\R,K)_\phi) 
\to \pi_0(K)^\phi \to\1. $$

Now we turn to the connecting maps. For that we note that 
for continuous sections, the map 
\begin{align*}
\Phi \: \Omega K &:= C_*(\bS^1,K) := C_*(\R/\Z,K)   
\to C_*(\R,K)_\phi, \\ 
\quad \Phi(f)(t) &:= \phi^{-n}(f([t])) \quad\mbox{ for } \quad 
t \in [n,n+1], n \in \Z 
\end{align*}
defines an isomorphism of Lie groups. 

\begin{proposition} For $j \geq 1$, the connecting maps 
\[
 \delta_j \: \pi_j(K) \to 
\pi_{j-1}(C_*(\R,K)_\phi) 
\cong \pi_{j-1}(\Omega K) \cong \pi_j(K), 
\]
are group homomorphisms given by 
$\delta_j([f]) = [f] - [\phi^{-1} \circ f].$
\end{proposition}

\begin{proof} For the adjoint action of $K$ on itself, this formula
is the Samelson product with $[k]\in\pi_{0}(K)\cong
\Bun(\bS^{1},K)$, and the proof in \cite[Thm.~2.4]{Wo07b} implies 
the present assertion when applied to $K \rtimes_\phi \Z$ instead of~$K$. 
\end{proof}

\begin{remark} \label{rem:3.6} (a) We have a short exact sequence 
$$\1\to \pi_1(K)_\phi := \pi_1(K)/\im(\pi_1(\phi) -\id) 
\into \pi_0(C^\infty(\R,K)_\phi) \onto 
\pi_0(K)^\phi \to\1. $$
If $K$ is connected, we obtain in particular 
$\pi_0(C^\infty(\R,K)_\phi) \cong \pi_1(K)_\phi.$

(b) For the evaluation of period maps, important information is contained 
in the short exact sequence 
$$ \1 \to \pi_3(K)_\phi \into \pi_2(C^\infty(\R,K)_\phi) \onto 
\pi_2(K)^\phi \to\1. $$
If $\pi_2(K)$ vanishes, it follows that the corresponding map 
$$ \pi_3(K) \cong \pi_2(\Omega K) \to \pi_2(C^\infty(\R,K)_\phi) $$
is surjective. 
\end{remark}

\begin{example} We discuss some examples where $\phi$ acts non-trivially 
on $\pi_2(K)$. 

(a) A typical example of a Lie group $K$ for which 
$\pi_2(K)$ is non-trivial is the projective unitary group 
$\PU({\cal H})$ of an infinite-dimensional 
complex Hilbert space ${\cal H}$ (\cite{Ku65}). 
Each automorphism of this simply 
connected group either is induced by a unitary or an anti-unitary map. 
In fact, the simple connectedness of $\PU({\cal H})$ implies 
$ \Aut(\PU({\cal H})) \cong \Aut(\pu({\cal H})),$
and since the Lie algebra $\fu({\cal H})$ is the universal 
central extension of $\pu({\cal H})$ 
(\cite[Example~III.6]{Ne02b}), each automorphism of 
$\pu({\cal H})$ lifts to a unique automorphism of 
$\fu({\cal H})$, so that 
$$\Aut(\pu({\cal H})) 
\cong \Aut(\fu({\cal H}))
\cong \Aut_*(\gl({\cal H})) 
\cong \PU({\cal H}) \rtimes \Z/2, $$
where the latter isomorphism follows from 
Prop.~3 in Section II.13 of \cite{dlH72} and 
$\Aut_*$ denote the group of all automorphism $\phi$ with 
$\phi(x^*) = \phi(x)^*$. 
We conclude that 
$$ \pi_0(\Aut(\PU({\cal H}))) \cong \Z/2. $$
Conjugation with an anti-unitary map induces the inversion on 
the center $\T \id_{\cal H}$ 
of $\UU({\cal H})$, and this implies that 
the action of $\pi_0(\Aut(\PU({\cal H})))$ 
induces the inversion on 
$\pi_2(\PU({\cal H})) \cong \pi_1(Z(U({\cal H})))\cong \Z$. 

(b) Another example is the smooth 
loop group $C^\infty(\bS^1,C)$ 
of a compact simple simply connected Lie group $C$ which satisfies 
$$ \pi_2(C^\infty(\bS^1,C)) \cong \pi_3(C) \cong \Z. $$
Its automorphism group is 
$$ \Aut(C^\infty(\bS^1,C))\cong C^\infty(\bS^1,\Aut(C)) \rtimes \Diff(\bS^1) $$
(cf.\ \cite[Prop.~3.4.2]{PS86}) 
whose group of connected components is 
$$ \big(\pi_1(\Aut(C)) \rtimes \pi_0(\Aut(C))\big) 
\rtimes \Z/2 \cong \big(Z(C) \rtimes \pi_0(\Aut(C))\big)\rtimes \Z/2, $$
hence finite.
In particular, any orientation reversing diffeomorphism acts on
$\pi_{2}(C^{\infty}(\bS^{1},C))\cong \Z$ by inversion.
\end{example}

\subsection{Period maps for twisted loop groups} \label{subsec:2.2}

Let $\kappa \: \fk \times \fk \to V$ be a $\fk$-invariant 
symmetric bilinear form and $\phi_V \in \GL(V)$ (defining a 
$\Z$-module structure on $V$) with 
$$ \phi_V(\kappa(x,y)) = \kappa(\phi_\fk x, \phi_\fk y) \quad \mbox{ for } 
\quad x,y \in \fk. $$
We write $\V := \V_{\phi_V}$ for the vector bundle over $\bS^1$ 
with fiber $V$ and holonomy $\phi_V$. 
Then the cocycle
corresponding to the canonical connection $\nabla$ on ${\cal K}_\phi$ 
defined by $\dd^\nabla f = f'(t)dt$ is given by 
$$ \omega_{\varphi }(f,g) 
:= [\kappa(f,\dd^\nabla g)] 
= \Big[\int_0^1 \kappa(f,g')\, dt\Big] \in 
V_{\phi_V} \cong \oline\Omega^1(\bS^1,\V), $$ 
where the last isomorphism comes from the following lemma:

\begin{lemma} \label{lem:a.5} Let $\V$ be a vector bundle over $\bS^1$ with fiber $V$ and 
holonomy $\phi_V \in \GL(V)$. Identifying $\bS^1$ with $\R/\Z$, the map 
$$ \oline\Omega^1(\bS^{1},\V) \cong H^1_{\rm dR}(\bS^1,\V) \to V_{\phi_V} 
= \coker(\phi_V - \id_V), \quad 
[f\cdot dt] \mapsto \Big[\int_0^1 f(t)\, dt\Big]$$  
for $f\in\Gamma \V \cong C^{\infty}(\R,V)_{\varphi_{V}}$ 
is a linear isomorphism. 
\end{lemma}

\begin{proof} Write 
$\frac{d}{dt}$ for the vector field generating the rigid rotations of 
$\bS^1$. 
Then ${\cal V}(\bS^1) = C^\infty(\bS^1,\R)\frac{d}{dt}$ 
implies that evaluation in $\frac{d}{dt}$ leads to an isomorphism 
$$\Omega^1(\bS^1,\V) 
= \Hom_{C^\infty(\bS^1,\R)}({\cal V}(\bS^1), \Gamma\V)
\cong \Gamma\V \cong C^\infty(\R,V)_{\phi_V}, $$
and under this identification, the 
canonical covariant derivative is given by 
$$ \dd^\nabla \: \Gamma \V \to \Omega^1(\bS^1,\V), \quad 
\dd^\nabla f = f'. $$
We first observe that if 
$f = g'$ for some $g \in \Gamma\V$, then 
$$ \int_0^1 f(t)\, dt = g(1) - g(0) = \phi_V^{-1}(g(0)) - g(0) 
\in \im(\phi_V^{-1} - \id_V) = \im(\phi_V - \id_V), $$
and the map is well-defined.
If, conversely, 
$f \: \R \to V$ is a smooth function, representing an element of 
$\Omega^1(\bS^1,\V)$, for which $\int_0^1 f(t)\, dt = \phi_V^{-1}(v)-v$ for some 
$v \in V$, then 
$$ g(t) := v + \int_0^t f(\tau)\, d\tau $$
satisfies $g' = f$ and 
\begin{align*}
g(t + 1) 
&= v + \int_0^1 f(\tau)\, d\tau + \int_1^{t+1} f(\tau)\, d\tau \\
&= \phi_V^{-1}(v) + \phi_V^{-1}\big(\int_0^{t} f(\tau)\, d\tau\big) 
= \phi_V^{-1}(g(t)). 
\end{align*}
This proves injectivity. To obtain surjectivity, choose some 
$\gamma \from
[0,1]\to[0,1]$ which is smooth, constant on a neighborhood of the boundary
and satisfies $\gamma (0)=0$, $\gamma (1)=1$ and $\int_{0}^{1}\gamma
(t)dt=\frac{1}{2}$. Then, 
for each $v\in V$, the mapping  
\mbox{$t\mapsto(1-\gamma (t))\cdot v+\gamma
(t)\varphi_{V}^{-1}(v)$} can be extended to an element $f_{v}$ of
$\Gamma \V$ with $[\int_{0}^{1}f_{v}(t)\, dt]=[v]$.
\end{proof}

\begin{remark} \label{rem:t2} If $V$ is infinite-dimensional, we further 
assume that the image of the operator $\phi_V - \id_V$ is closed, so that 
$\oline\Omega^1(\bS^1,\V_{\phi_V}) \cong V_{\phi_V}$ is Hausdorff. 
\end{remark}

To study the period map (cf.\ Appendix C)\begin{footnote}{Note 
that we do not have to impose any completeness 
condition on the quotient space $V_{\phi_V}$ 
to make sense of the period 
integrals because they can be calculated as $V$-valued integrals.}
\end{footnote}   
$$ \per_{\omega_{\varphi }}\: \pi_2(C^\infty(\R,K)_\phi) \to 
V_{\phi_V}, $$
we first consider the subgroup 
$C^\infty_*(\R,K)_\phi$. 
Since $C^\infty(\R,K)_\phi$ is locally exponential (Appendix A), 
this is a Lie subgroup with Lie algebra 
$C^\infty_*(\R,\fk)_\phi :=\ker \ev_{0}^\fk$. 
To evaluate the period map on 
$\pi_2(C^\infty_*(\R,K)_\phi)$, we note that on 
$C^\infty_*(\R,\fk)_{\varphi },$
$$ \tilde\omega_{\varphi }(f,g) := \int_0^1 \kappa(f,g')(t)\, dt \in V $$
defines a Lie algebra cocycle. In fact, integration by parts shows that 
it is alternating, and with Remark~\ref{rem:1.2} we obtain 
for $f,g,h \in C^\infty_*(\R,\fk)_{\varphi }$
$$ -(\dd\tilde\omega_{\varphi })(f,g,h) 
= \int_0^1 \kappa([f,g],h)'
= \kappa([f,g],h)(1) - \kappa([f,g],h)(0) =0. $$

The following lemma reduces the period map of $\tilde\omega_\phi$ 
(cf.\ Appendix \ref{app:C}) to the more accessible period map 
$\per_{C(\kappa)^{l}}\from \pi_{3}(K)\to V$ 
of the closed 
$3$-form $C(\kappa)^l$ on $K$ which is studied in 
detail in Appendix~\ref{app:B} for $\dim K < \infty$. 

\begin{lemma} \label{lem:per-loop} 
Identifying $\pi_3(K)$ in the canonical way with the group $\pi_2(\Omega K) \cong 
\pi_2(C^\infty_*(\R,K)_\phi)$, we have 
$$ \per_{\tilde\omega_{\varphi }} = \frac{1}{2}\per_{C(\kappa)}\: 
\pi_3(K) \to V. $$
More generally, if $\sigma \: \bS^2\to C^\infty(\R,K)_\phi$ is a smooth 
map and $\tilde\sigma \: \R \times \bS^2 \to K$ defined by 
$\tilde\sigma(t,m) := \sigma(m)(t)$, then 
\begin{equation}
  \label{eq:per-3form}
\per_{\omega_\phi}([\sigma]) = \frac{1}{2} \Big[\int_{[0,1] \times \bS^2} 
\tilde\sigma^*C(\kappa)^l\Big]. 
\end{equation}
\end{lemma}

\begin{proof} It suffices to verify that the $V$-valued 
Lie algebra $2$-cochain \break 
$\tilde\omega_\phi(f,g) := \frac{1}{2} \int_0^1 \kappa(f,g') - \kappa(g,f')\, dt$ satisfies  
$$ \int_\sigma \tilde\omega_\phi^l = \frac{1}{2} 
\int_{[0,1] \times \bS^2} \tilde\sigma^*C(\kappa)^l.$$
Since homotopy classes may be represented by smooth maps 
\cite[Sect.~A.3]{Ne02a}, both assertions follow from that.

First we note that $\sigma$ also defines a smooth curve 
in $\hat \sigma \in 
C^\infty(\R, C^\infty(\bS^2,K))$ by 
$\hat\sigma(t)(m) := \sigma(m)(t)$. We then identify 
its logarithmic derivative 
$\delta(\hat\sigma)$ with a smooth curve with values in 
$\L(C^\infty(\bS^2,K)) = C^\infty(\bS^2,\fk)$, 
so that $\dd \delta(\hat\sigma)$ is a smooth curve with values in 
$\Omega^1(\bS^2,\fk)$. 

We consider 
$\delta\sigma \in \Omega^1(\bS^2, C^\infty(\R,\fk)_\phi) 
\cong C^\infty(\R, \Omega^1(\bS^2, \fk))$ 
as a smooth curve with values in 
$\Omega^1(\bS^2,\fk)$ in the obvious fashion. 
Using the fact that $\delta \tilde\sigma \in \Omega^1(\R \times \bS^2,\fk)$ 
satisfies the Maurer--Cartan equation 
$$ \dd \delta\tilde\sigma 
+ \frac{1}{2}[\delta\tilde\sigma, \delta\tilde\sigma] = 0, $$
the derivative of this curve can be 
calculated by evaluating it on some smooth vector field $X \in \cV(\bS^2)$: 
\begin{align*}
(\delta\sigma)'(X) 
&= {\cal L}_{\frac{\partial}{\partial t}}(\delta\tilde\sigma(X)) 
= \big({\cal L}_{\frac{\partial}{\partial t}}\delta\tilde\sigma\big)(X) 
= \Big(i_{\frac{\partial}{\partial t}}\dd\delta\tilde\sigma
+ \dd i_{\frac{\partial}{\partial t}}\delta\tilde\sigma\Big)(X) \\ 
&= \Big(-\frac{1}{2}i_{\frac{\partial}{\partial t}}[\delta\tilde\sigma, 
\delta\tilde\sigma]+ \dd \delta(\hat\sigma)\Big)(X)
= -[\delta\tilde\sigma(\frac{\partial}{\partial t}), 
\delta\tilde\sigma(X)]+ \dd \delta(\hat\sigma)(X) \\
&= -[\delta\hat\sigma, \delta\sigma(X)]+ \dd \delta(\hat\sigma)(X) 
= \big(\dd \delta(\hat\sigma)+ [\delta\sigma, \delta\hat\sigma]\big)(X). 
\end{align*}
This proves that 
\begin{equation}
  \label{eq:dif-rel}
(\delta\sigma)' = \dd \delta(\hat\sigma)+ [\delta\sigma, \delta\hat\sigma] 
\in C^\infty(\R,\Omega^1(\bS^2,\fk)), 
\end{equation}
Using the Maurer--Cartan equation for $\delta\sigma$, we further get 
in $C^\infty(\R, \Omega^2(\bS^1,\fk))$: 
\begin{align*}
\delta \sigma \wedge_\kappa (\delta\sigma)' 
&= \delta \sigma \wedge_\kappa \dd \delta(\hat\sigma) 
+ \delta\sigma \wedge_\kappa [\delta\sigma, \delta\hat\sigma] 
= \delta \sigma \wedge_\kappa \dd \delta(\hat\sigma) 
+ [\delta\sigma, \delta\sigma]  \wedge_\kappa \delta\hat\sigma \\
&= - \dd(\delta\sigma \wedge_\kappa \delta\hat\sigma\big)
+  \dd(\delta\sigma)\wedge_\kappa \delta\hat\sigma 
+ [\delta\sigma, \delta\sigma]  \wedge_\kappa \delta\hat\sigma \\
&= - \dd(\delta\sigma \wedge_\kappa \delta\hat\sigma\big)
- \frac{1}{2} [\delta\sigma, \delta\sigma] \wedge_\kappa \delta\hat\sigma 
+ [\delta\sigma, \delta\sigma]  \wedge_\kappa \delta\hat\sigma \\
&= - \dd(\delta\sigma \wedge_\kappa \delta\hat\sigma\big)
+ \frac{1}{2} [\delta\sigma, \delta\sigma] \wedge_\kappa \delta\hat\sigma\\ 
&= - \dd(\delta\sigma \wedge_\kappa \delta\hat\sigma\big)
+ C(\kappa)(\delta\sigma, \delta\sigma,\delta\hat\sigma)\\ 
\end{align*}
and hence  
\begin{align*}
\int_\sigma \tilde\omega_\phi^l 
&= \int_{\bS^2} \tilde\omega_\phi(\delta\sigma,\delta\sigma)
= \frac{1}{2}\int_0^1 \int_{\bS^2} \delta\sigma \wedge_\kappa (\delta\sigma)'\, dt \\ 
&= \frac{1}{2}\int_0^1\int_{\bS^2}  
C(\kappa)(\delta\sigma,\delta\sigma,\delta\hat\sigma) \, dt
= \frac{1}{2}\int_{[0,1] \times \bS^2} \tilde\sigma^*C(\kappa)^l.
\\\end{align*}
For the last equality we have used that 
$\tilde\sigma^*C(\kappa)^l = C(\kappa)(\delta \sigma, 
\delta\sigma, \delta\hat\sigma)\, dt$, 
which is most easily verified by applying 
both sides to triples of tangent vectors of the form 
$(\frac{\partial}{\partial t},v,w)$ for $v,w \in T_m(\bS^2)$. 
\end{proof}

The preceding lemma shows in particular that the period homomorphism 
$\per_{\tilde\omega_\phi}$ 
does not depend on the pair $(\varphi, \varphi_V)$. 

\begin{example} \label{ex:su2} It is instructive to take a closer look at the 
example $K = \SU_2(\C)$. We realize $\SU_2(\C) \cong \bS^3$ 
as the group of unit quaternions in 
$\bH$ and write $\kappa(x,y) = -\frac{1}{4}\tr(\ad x\ad y)$ for the normalized 
invariant symmetric bilinear form, satisfying 
$\kappa(x,x) = 2 \|x\|^2$ for each $x \in \su_2(\C) = \Spann_\R\{I,J,K\}$. 
For the basis elements $I,J,K$, we then have 
$\kappa([I,J],K) = 2\kappa(K,K) = 4,$
so that the left invariant $3$-form defined by $C(\kappa)(x,y,z) := 
\kappa([x,y],z)$ on $\SU_2(\C) \cong \bS^3$ is $4 \mu_{\bS^3}$, where 
$\mu_{\bS^3}$ is the volume form of $\bS^3$. It follows in particular, that 
\begin{equation}
  \label{eq:3-period}
\per_{C(\kappa)^l}([\id_K]) 
= \int_{\SU_2(\C)} C(\kappa)^l = 4 \vol(\bS^3) = 8 \pi^2.
\end{equation}
On the other hand, it has been shown in \cite{PS86} 
(see also the calculations in Appendix IIa to Section II in \cite{Ne01}) 
that 
$\frac{1}{2\pi} \Pi_{\omega_\kappa} 
= \Pi_{\frac{1}{2\pi}\omega_\kappa} = 2\pi \Z$
for $\omega_{\kappa}:=\wt{\omega}_{\id}$,
so that 
$$ \Pi_{\omega_\kappa} = 4\pi^2 \Z. $$
In view of the preceding lemma, this is a direct consequence of 
\eqref{eq:3-period}. 
\end{example}

Let $q_V \: V \to V_{\phi_V}$ denote the projection map. 
Then the cocycle $\tilde\omega_\phi$ satisfies 
$$ q_V \circ \tilde\omega_\phi = \L(\iota)^*\omega_\phi, $$
where $\iota \: C^\infty_*(\R,K)_\phi \into 
C^\infty(\R,K)_\phi$ denotes the inclusion map. 
Therefore Remark~\ref{rem:8.4} yields 
\begin{equation}
  \label{eq:pre-rel3}
\per_{\omega_\phi} \circ \pi_2(\iota) 
= \per_{\L(\iota)^*\omega_\phi} 
= q_V \circ \per_{\tilde\omega_\phi} 
= q_V \circ \per_{C(\kappa)}. 
\end{equation}
If $\pi_2(K)$ vanishes, then $\pi_2(\iota)$ is surjective 
(Remark~\ref{rem:3.6}) and we thus obtain: 

\begin{theorem} If $\pi_2(K)$ vanishes, then 
$\Pi_{\omega_\phi} \subeq q_V\big(\im(\per_{C(\kappa)})\big).$
\end{theorem}

As a consequence, we obtain for finite dimensional groups with 
Theorem~\ref{thm:disc1} and Cartan's Theorem that 
$\pi_2(K)$ vanishes in this case (Remark~\ref{rem:mn.2.3}). 

\begin{corollary} \label{cor:disc-univ} If $K$ is finite-dimensional, 
$V= V(\fk)$ and $\kappa = \kappa_u$ is universal, then the period group 
$\Pi_{\omega_\phi}$ is discrete. 
\end{corollary}

We now present an example where the period group depends 
significantly on the connection $\nabla$. 

\begin{example} \label{ex:torus} 
 We consider the special case of Remark~\ref{rem:princbun-ex}, where  
$$ \pi \: Q = \T^2 \to M = \T = \R/\Z, \quad \pi(t,s) = s, $$
$H = \T$ and $K = C^\infty(\T,G)$ for a compact simple Lie group~$G$. 
Then $\Gamma\cK \cong C^\infty(\T^2,G)$ and 
$\cK \cong \T \times K$ is a trivial bundle. 

To a positive definite invariant symmetric bilinear form 
$\kappa_\g$ on $\g$, we associate the invariant bilinear form 
$$ \kappa(f,g) := \int_0^1 \kappa_\g(f(t),g(t))\, dt $$ 
on $\fk$. 
The group $H = \T$ acts on $K$, resp., $\fk$, by composition, 
and the action of the Lie algebra $\fh \cong \R$ is given by 
$Df = f'$, which leaves $\kappa$ invariant. The
Lie algebra cocycle
$$\eta_D(f,g) := \kappa(f,Dg) 
= \int_0^1 \kappa_\g(f(t), g'(t))\, dt, $$
on $\fk$ is universal (cf.\ \cite{PS86}). 
In particular, the period 
homomorphism 
$$ \per_{\eta_D} = \frac{1}{2}\per_{C(\kappa_\g)} \: 
\pi_2(K) \cong \pi_3(G) \cong \Z \to \R $$
is non-trivial. 

The covariant exterior derivatives on the trivial bundle 
$P = \bS^1 \times H$ take on $\Gamma\fK \cong C^\infty(\T,\fk) 
\cong C^\infty(\T^2,\g)$ the form 
$$ \dd^\nabla f = f' + h\cdot Df 
= \frac{\partial f}{\partial s} + h(s)\cdot \frac{\partial f}{\partial t}, $$
for some $h\in C^{\infty }(\T,\R)$, determined by $\nabla$.
Accordingly, the cocycle $\omega_\kappa^\nabla$ decomposes as 
$$ \omega_\kappa^\nabla(f,g) = 
\underbrace{\int_0^1 \kappa(f,\pds{g})\, ds}_{\omega_0(f,g):=}
 + \underbrace{\int_0^1 h(s)\kappa(f,\pdt{g})\, ds}_{\eta(f,g):=}.  
$$

To calculate the period maps for the cocycles
$$ \omega_0(f,g) = \int_0^1\int_0^1 \kappa_\g(f,\pds{g})\, ds\, dt, \quad 
\eta(f,g) = \int_0^1\int_0^1 h(s)\kappa_\g(f,\pdt{g})\, ds\, dt, $$
we express them in terms of the universal cocycle 
$$ \omega_u(f,g) = [\kappa_\g(f,\dd g)] 
= [ \kappa_\g(f, \pdt g)\, \dd t + \kappa_\g(f, \pds g)\, \dd s] $$
of $\Gamma\fK 
= C^\infty(\T^2,\g)$ with values in $\oline\Omega^1(\T^2,\R)$. 
If $\kappa_\g$ is suitably normalized, the period group of this cocycle 
is $\Pi_{\omega_u} =\Z [\dd t] + \Z [\dd s]$ (\cite{MN03}). 
We further find with Remark~\ref{rem:mn.2.3}
$$ \pi_2(C^\infty(\T^2,G)) 
\cong \pi_2(G) \oplus \pi_3(G)^2 \oplus \pi_4(G)  
\cong \Z^2 \oplus \pi_4(G) $$ 
(cf.\ \cite[Rem.~I.11]{MN03}), and in these terms 
$$ \per_{\omega_u}(m,n,u) = m [\dd t] + n [\dd s].$$

If $\gamma_t(s) = \gamma^s(t) = (t,s)$ describes the vertical 
and horizontal circles in $\T^2$, then
$$ \omega_0(f,g) = \int_0^1 \cI_{\gamma_t} \circ \omega_u(f,g)\, dt,  $$
where $\cI_{\gamma_t}[\alpha] := \int_{\bS^1} \gamma_t^*\alpha$, 
so that the period map is given by 
$$ \per_{\omega_0}(m,n,u)
 = \int_0^1 \cI_{\gamma_t} \circ \per_{\omega_u}(m,n,u)\, dt 
 = \int_0^1 n\, dt = n. $$
Similarly, 
$$ \eta(f,g) = \int_0^1 h(s) \cI_{\gamma_s} \circ \omega_u(f,g)\, ds,  $$
and its period map is 
$$ \per_{\eta}(m,n,u)
 = \int_0^1 h(s) \cI_{\gamma_s} \per_{\omega_u}(m,n,u)\, ds 
 = \int_0^1 h(s)m\, ds = m \int_0^1 h(s)\, ds. $$
This implies that the period group 
$$\Pi_\omega = \Z + \Z\cdot \int_0^1 h(s)\, ds $$
is discrete if and only if the integral $\int_0^1 h(s)\, ds$ is rational. 
\end{example}

\begin{remark}
We have seen in Remark~\ref{rem:3.6} that for a twisted loop group 
$L_\phi K := C^\infty(\R,K)_\phi$, $\phi \in \Aut(K)$, 
the group $\pi_2(L_\phi K)$ is determined by a short exact sequence 
$$ \1 \to \pi_3(K)_\phi \to \pi_2(L_\phi K) \to \pi_2(K)^\phi \to \1.$$
Accordingly, the period group $\Pi_{\omega_\phi}$ can be determined 
in a two-step process. The restriction to $\pi_3(K)_\phi$ is, 
up to the factor $\frac{1}{2}$,  
the period map 
$$ \per_{C(\kappa),\phi} \: \pi_3(K)_\phi \to V_{\phi_V}, \quad 
[\sigma] \mapsto [\per_{C(\kappa)}(\sigma)] $$
obtained by factorization of $\per_{C(\kappa)}$. If 
$\Pi_{C(\kappa),\phi} \subeq V_{\phi_V}$ denotes the image of this 
homomorphism, then $\per_\omega$ factors through a homomorphism 
$$ \oline\per_\omega \: \pi_2(K)^\phi \to V_{\phi_V}/\Pi_{C(\kappa),\phi}$$
whose image determines the period group $\Pi_\omega$ as its inverse 
image in $V_{\phi_V}$. 
\end{remark}

The following example shows that both parts $\pi_3(K)_\phi$ and 
$\pi_2(K)^\phi$ may contribute non-trivially to $\Pi_\omega$ and that 
the period group depends seriously on $\phi$. 

\begin{example} (a) Let $G$ be a simply connected simple compact Lie group and 
$K := C^\infty(\bS^1,G)$ be its loop group.

Let 
$$\phi_K = (h,\psi)
\in \Aut(K) \cong C^\infty(\bS^1,\Aut(G)) \rtimes \Diff(\bS^1) $$
(cf.\ \cite[Prop.~3.4.2]{PS86}). 
Here $\Aut(K)$ actually carries a natural 
Lie group structure and the automorphism $\phi_V$ of 
$V =\R$ induced by $\phi_K$ for which the form 
$\kappa(f,g) = \int_{\bS^1} \kappa_\g(f(t), g(t))\, dt$ 
is invariant is $\pm \id_V$, 
depending on whether $\psi$ preserves the orientation 
of $\bS^1$ or not. 
We also note that 
$$\pi_0(C^\infty(\bS^1,\Aut G)) \cong \pi_1(\Aut G) \rtimes \pi_0(\Aut G)  
\cong Z(G) \rtimes \pi_0(\Aut G) $$
is a finite group and that the subgroup 
$$\pi_0(C^\infty_*(\bS^1,\Aut G)) \cong \pi_1(\Aut G) $$ 
acts trivially on all higher homotopy groups of 
$G$,\begin{footnote}{Here we use that if a topological group acts 
on a space $M$, then the corresponding action of $\pi_1(G)$ on 
$\pi_k(M,x_0)$, $k \geq 1$, is always trivial. 
One finds the special cases where $G$ acts on 
itself by the multiplication map in \cite{Hu59}, Prop.~16.10. 
The general case is proved similarly.}
\end{footnote}hence in particular 
on $\pi_3(G) \cong \pi_2(K)$. Moreover, $\Aut(G)$ preserves the 
Cartan--Killing form $\kappa_\g$ of $\g$, hence fixes the associated closed 
invariant $3$-form, so that de Rham's Theorem 
implies that it also acts trivially on $\pi_3(G)$. 

For $\phi_V = -\id_V$ we obtain in particular 
$V_{\phi_V} = \{0\}$, so that all periods vanish. In the latter 
case, the natural identification of $\pi_2(K)$ with $\pi_3(G)$ 
shows that $\phi_K$ acts as $-\id$ on 
$\pi_2(K) \cong \pi_3(G) \cong \Z$, so that $\pi_2(K)^\phi = \{0\}$. 

If $\phi_V = \id_V$, then $V_{\phi_V} = V = \R$ and $\psi$ is orientation 
preserving. Then the action of $\phi_K$ on $\pi_2(K)\cong \pi_3(G)\cong \Z$ 
is trivial. The action of $\phi_K$ on 
$\pi_3(K) \cong \pi_3(G) \oplus \pi_4(G)$ is trivial on the first factor 
(coming from constant functions) and $\pi_4(G)$ is finite, so that 
$\pi_3(K)_{\phi}$ is of rank~$1$. 
Now the same arguments as in 
Example~\ref{ex:torus} show that $\pi_2(L_\phi K)$ is of 
rank $2$ and both summands contribute to~$\Pi_\omega$. 
\end{example}

\section{Corresponding Lie group extensions} \label{sec:3}

We now determine in which cases the central extension $\wh{\Gamma
\fK}$ defined by the cocycle $\omega_{\kappa}^{\nabla }$ integrates to
a Lie group extension. To this end we analyze its period group
$\Pi_{\omega}:=\im(\per_{\omega_{\kappa}^{\nabla}})$ and determine
whether the adjoint action of $\Gamma \cK$ lifts to an action on
$\wh{\Gamma \fK}$ (cf.\ Appendix~\ref{app:C} and \cite{Ne02a}). 
Throughout, $M$ denotes a compact connected manifold. 

\subsection{On the image of the period map}

Throughout this section, we fix a base point 
$p_0 \in P$ and put $m_0 := q_P(p_0)$. We also assume that the Lie group $H$ 
is regular. 

It is convenient to consider an intermediate 
situation given by a covering manifold $\hat q_M \: \hat M \to M$,  
defined as follows. 
Let $\delta_1 \: \pi_1(M) \to \pi_0(H)$ denote the connecting map 
from the long exact homotopy sequence of the 
principal $H$-bundle $P$ that we used to define ${\cal K}$. We write 
$$\tilde\rho_V := \rho_V \circ \delta_1\: \pi_1(M) \to\GL(V)$$ 
for the corresponding pullback representation of $\pi_1(M)$ on $V$ 
and put $\hat M := \tilde M/\ker \tilde\rho_V$. Then 
$\hat M$ is a covering of $M$ with $\pi_1(\hat M) \cong \ker \tilde\rho_V$, 
and its group of deck transformations is $\D := \pi_1(M)/\pi_1(\hat M) 
\cong \tilde\rho_V(\pi_1(M))$. 

Since $P$ is connected, the connecting homomorphism $\delta_1$ is 
surjective and the squeezed bundle $P/H_0$ is a covering of 
$M$ associated to $\delta_1$, hence equivalent to 
$\tilde M \times_{\delta_1} \pi_0(H) \cong \tilde M/\ker \delta_1$. 
This implies that 
$$P/\ker\rho_V \cong (P/H_0)/\ker\oline\rho_V 
\cong \tilde M/\ker \tilde\rho_V \cong \hat M. $$

\begin{remark} For the open subgroup $H_V := \ker \rho_V$ of $H$, 
we may also consider $P$ as an principal $H_V$-bundle 
$\hat q \: P \to \hat M$. 
If $\hat\V := \hat q_M^*\V$ denotes the pullback of $\V$ to $\hat M$, 
it follows that 
$\hat\V \cong P \times_{\rho_V\res_{H_V}} V \cong \hat M \times V$ 
is a trivial vector bundle, which leads to a natural map 
$$ \oline\Omega^1(M,\V) \to \oline\Omega^1(\hat M,V). $$
\end{remark}

In the following our first step to the understanding of the period group 
of $\omega := \omega_\kappa^\nabla$ is to investigate when 
it is contained in the subspace 
$H^1_{\rm dR}(M,\V)$ of $\oline\Omega^1(M,\V)$. If this condition is 
not satisfied, then one may not expect any simple criteria for 
discreteness, as the  Examples~\ref{ex:counter-ex-h1} and \ref{ex:torus} 
show. 

\begin{definition} \label{def:6.1} 
\nin (a) Fix a connection $\nabla$ on the principal $H$-bundle~$P$. 
For any smooth loop $\alpha \: [0,1] \to M$, based in 
$m_0 \in M$, we define its {\it holonomy} ${\cal H}_{p_0}(\alpha) \in H$ as 
follows. Since $H$ is assumed to be regular, the curve 
$\alpha$ has a unique smooth horizontal lift 
$\hat\alpha \: [0,1] \to P$ starting in $p_0$ 
(cf.\ \cite{KM97}), and since 
$\hat\alpha(1)$ and $\hat\alpha(0)$ are both mapped to 
$\alpha(0) = \alpha(1) = m_0$, there exists a unique element 
${\cal H}_{p_0}(\alpha) \in H$ with 
$$ \hat\alpha(1) = \hat\alpha(0).{\cal H}_{p_0}(\alpha). $$
Changing the base point leads to the relation
$$ {\cal H}_{p_0.h}(\alpha) = h^{-1} {\cal H}_{p_0}(\alpha)h,$$  
so that the holonomy depends on $p_0$. Since we keep the base point 
$p_0$ fixed, we may also write ${\cal H}(\alpha) := {\cal H}_{p_0}(\alpha)$. 
If $\alpha \: \R \to M$ is a $1$-periodic map with 
$\alpha(0)= m_0$, representing a smooth loop 
$\R/\Z \to M$, then we put ${\cal H}(\alpha) 
:= {\cal H}(\alpha\res_{[0,1]})$.

\nin (b) We identify the group $\Gamma \cK$ and the Lie algebra 
$\Gamma \fK$ with the corresponding 
spaces of $H$-equivariant maps $C^{\infty}(P,K)^{H}$, resp., 
$C^{\infty}(P,\fk)^{H}$. 

Any smooth $1$-periodic map $\alpha \: \R \to M$ lifts to a unique 
smooth horizontal curve $\hat \alpha \: \R \to P$ satisfying 
$$ \hat\alpha(t+1)  = \hat\alpha(t).{\cal H}(\alpha) $$
for each $t \in \R$ (both sides define horizontal curves and coincide 
in $t = 0$). We put 
$\phi_K^\alpha := \rho_K({\cal H}(\alpha))$ and 
$\phi_\fk^\alpha := \L(\phi_K^\alpha)$. Then we obtain a homomorphism of 
Lie groups 
\[
\hat\alpha_K^* \from \Gamma \cK \to C^\infty(\R,K)_{\phi_K^\alpha}, 
\quad f \mapsto f \circ \hat\alpha, 
\]
and of Lie algebras 
\[
\hat\alpha_\fk^* = \L(\hat\alpha_K^*)
 \from \Gamma \fK \to C^\infty(\R,\fk)_{\phi_\fk^\alpha}, 
\quad f \mapsto f \circ \hat\alpha. 
\]

\nin (c) As before, we realize $\Omega^1(M,\V)$ as the 
space $\Omega^1(P,V)^H$ of $H$-invariant basic $V$-valued $1$-forms on $P$. 
For each smooth loop $\alpha$ in $m_0$ we put 
$\phi_V^\alpha := \rho_V({\cal H}(\alpha))$. For each 
$H$-equivariant smooth function $f \: P \to V$ we then have 
\begin{align*}
\int_0^1 \hat\alpha^*\dd f 
&= f(\hat\alpha(1)) - f(\hat\alpha(0)) 
= {\cal H}(\alpha)^{-1}.f(\hat\alpha(0)) - f(\hat\alpha(0)) \\ 
&\in \im((\phi_V^\alpha)^{-1} - \id_V)
= \im(\phi_V^\alpha - \id_V). 
\end{align*}
Therefore we have a well-defined integration map 
$$ \cI_\alpha \: \oline\Omega^1(M,\V) \to V_{\phi_V^\alpha}
=\coker(\phi_V^\alpha-\id_{V}), 
\quad [\theta] \mapsto \Big[\int_0^1 \hat\alpha^*\theta\Big]. $$

\nin (d) Let $\hat q \: P \to \hat M$ be the bundle projection and 
$\hat\alpha_M := \hat q \circ \hat\alpha$. This is a piecewise smooth 
(continuous) lift of $\alpha$ to the covering space $\hat M$, starting in the 
base point $\hat m_0 := \hat q(p_0)$. Since the fibers of $\hat q$ 
are the orbits of $\ker \rho_V$, the condition $\phi_V^\alpha = \id_V$ 
is equivalent to the path $\hat\alpha_M$ being closed. If this is the 
case, then $\cI_\alpha$ has values in $V_{\phi_V^\alpha} = V$. 
\end{definition}

\begin{remark} Let $\omega_{\phi_K^\alpha}$ denote the canonical 
$V_{\phi_V}$-valued cocycle 
on $C^\infty(\R,\fk)_{\phi_\fk^\alpha}$ (cf.\ Subsection~\ref{subsec:2.2}). 
For the horizontal lift $\hat\alpha \: \R \to P$ and 
$f \in \Gamma\fK \cong C^\infty(P,\fk)^H$, we have 
$$ (\dd^\nabla f)(\hat\alpha'(t)) 
= \dd f(\hat\alpha'(t)) = (f \circ \hat\alpha)'(t). $$
Therefore 
\begin{align*}
\big(\L(\hat\alpha_K^*)^*\omega_{\phi_K^\alpha}\big)(f,g) 
&= \Big[ \int_0^1 \kappa(f \circ \hat\alpha, (g \circ \hat\alpha)')\, dt\Big]
= \Big[ \int_0^1 \hat\alpha^*\kappa(f,\dd^\nabla g)\Big]\\
&= \cI_\alpha([\kappa(f,\dd^\nabla g)]).
\end{align*}
We thus obtain the important relation 
\begin{equation}
  \label{eq:pull-rel}
\L(\hat\alpha_K^*)^*\omega_{\phi_K^\alpha}
= \cI_\alpha \circ \omega, 
\end{equation}
and with Remark~\ref{rem:8.4}, this leads to 
\begin{equation}
  \label{eq:pull-rel2}
\cI_\alpha \circ \per_{\omega} 
= \per_{\cI_\alpha \circ \omega} 
= \per_{\omega_{\phi_K^\alpha}} \circ \pi_2(\hat\alpha_K^*). 
\end{equation}
\end{remark}

\begin{remark} \label{rem:4.3} Let $F\from [0,1]\times \bS^{1} \to M$ be 
a smooth map which is a homotopy of loops based in $m_0$. 
Then 
$$\beta \from [0,1]\to H, \quad \beta(t) := {\cal H}(F_0)^{-1}
{\cal H}(F_{t})$$ 
is a smooth curve starting in $\1$. 
If $F$ is chosen to be independent of the first variable on a neighborhood of $\{0,1\}\times\bS^{1}$, then $\beta$ is constant on a neighborhood of $\{0,1\}$, 
so that it can be extended to a smooth map 
$$\beta \from \R\to H \quad \mbox{ with } \quad 
\beta (t+1)={\cal H}(F_0)^{-1}\beta (t){\cal H}(F_1) \quad \mbox{ for }\quad 
t \in \R. $$
For $i = 0,1$, put $\phi_i := \rho_K({\cal H}(F_i))$. 
Then 
\[
\Phi_{\beta} \from C^\infty(\R,K)_{\varphi_{1}}\to 
C^\infty(\R,K)_{\varphi_{0}},\quad
\Phi_{\beta} (f)(t)= (\beta.f)(t) := \rho_K(\beta(t))(f(t))
\]
is an isomorphism of Lie groups. The corresponding isomorphism
on the Lie algebra level is similarly given by 
$\big(\L(\Phi_{\beta})\xi\big)(t) = \rho_\fk(\beta(t)).\xi(t)$. 
The curve $\phi_{V,t} 
:= \rho_V({\cal H}(F_t))$ in $\GL(V)$ is constant 
because $H_0$ acts trivially. We may thus put $\phi_V := \phi_{V,t}$, 
and the target spaces of the cocycles 
$$\omega_{\phi_i}(f,g) = \Big[\int_0^1 \kappa(f,g')\, dt\Big] \in 
V_{\phi_V} $$
on $C^\infty(\R,K)_{\phi_i}$ coincide. 
Unfortunately, $\omega_{\phi_1}$ does not coincide with 
$\L(\Phi_\beta)^*\omega_{\phi_0}$. Instead, the product rule 
$(\beta g)' = \beta.(\delta(\beta).g) + \beta g',$ 
and the $H_{0}$-invariance of $\kappa$ show that 
$$ \big(\L(\Phi_\beta)^*\omega_{\phi_0}- \omega_{\phi_1}\big)(f,g) 
= \Big[\int_0^1 \kappa(f,\delta(\beta).g)\, dt\Big]. $$
Here we identify $\Omega^1(\R,\fh)$ with $C^\infty(\R,\fh)$, 
so that $\delta(\beta)$ is interpreted as a smooth $\fh$-valued 
curve. 

(b) If $q_V \: V \to V_{\phi_V}$ denotes the projection map, 
this means that 
\begin{equation}
  \label{eq:pull-rel3}
\L(\Phi_\beta)^*\omega_{\phi_0}- \omega_{\phi_1}
= q_V \circ \eta_{\delta(\beta)}, 
\end{equation}
where we put 
$$ \eta_\gamma(f,g) := \int_0^1 \kappa(f,\gamma.g)\, dt
\quad \mbox{ for } \quad \gamma \in C^\infty([0,1], \fh). $$
Since $\fh$ preserves 
$\kappa$, each $\eta_\gamma$ is a $V$-valued $2$-cocycle; actually 
$$ \eta_\gamma = \int_0^1 \eta_{\gamma(t)}\ dt 
\quad \mbox{ for } \quad 
\eta_{\gamma(t)}(x,y) := \kappa(x,\gamma(t).y), \quad 
\eta_{\gamma(t)} \in Z^2(\fk,V).$$

(c) Let $0 < \eps < \frac{1}{2}$. In $C^\infty_*(\R,K)_{\phi_1}$ we write 
$C^\infty_\eps(\R,K)_{\phi_1}$ for the subgroup of those maps 
vanishing on the interval $[-\eps,\eps]$. 
From Corollary~\ref{cor:3.4}, it easily follows that the inclusion 
of $C^\infty_\eps(\R,K)_{\phi_1}$ into
$C^{\infty}_{*}(\R,K)_{\varphi_{1}}$is a weak homotopy equivalence.
On the other hand, restriction to $[0,1]$ and periodic extension 
yields an isomorphism of Lie groups 
$$C^\infty_\eps(\R,K)_{\phi_1} \to C^\infty_\eps(\R,K)_{\id}. $$
Further, the isomorphism $\Phi(\beta)$ induces an automorphism 
of $C^\infty_\eps(\R,K)_{\id}$. With $\beta_t(s) := \beta(st)$, we even 
obtain by $\Phi(\beta_t)$ a smooth family of automorphisms of 
$C^\infty_\eps(\R,K)_{\id}$ connecting $\Phi(\beta)$ to the identity. 
Therefore $\Phi(\beta)$ induces the identity on 
$\pi_2\big(C^\infty_\eps(\R,K)_{\id}\big)$, and 
Lemma~\ref{lem:per-loop} implies that the period maps of 
$\tilde\omega_{\phi_1}$ and $\L(\Phi_\beta)^*\tilde\omega_{\phi_0}$ 
coincide. We conclude that the period map of the cocycle 
$\eta_{\delta(\beta)}$ vanishes on the image of 
$\pi_2(C^\infty_*(\R,K)_{\phi_1})$. If, in addition, $\pi_2(K)$ 
vanishes, then the long exact homotopy sequence of the map 
$\ev_0^K$ shows that the period map of 
$q_V \circ \eta_{\delta(\beta)}$ vanishes. 
\end{remark}

\begin{lemma} Let $\alpha_i \:\bS^1 \to M$, $i =0,1$, be two smooth homotopic 
loops in $m_0 \in M$ and $\beta \: [0,1] \to H$ a smooth curve in $H$ 
obtained from a smooth homotopy of $\alpha_0$ and $\alpha_1$ 
as in the preceding remark. Then the two morphisms of Lie groups 
$$ \hat\alpha_0^*,\quad \Phi_\beta \circ \hat\alpha_1^* \: 
\Gamma\cK \to C^\infty(\R,K)_{\phi_0} $$
are homotopic. 
\end{lemma}

\begin{proof} Let $F \: [0,1] \times \bS^1 \to M$ 
be a smooth homotopy with $F_i = \alpha_i$ for $i = 0,1$, 
and assume w.l.o.g.\ that $F$ is constant in a neighborhood 
of $\{0,1\} \times \bS^1$. Define $\beta$ as above. 
For $f  \in \Gamma\cK \cong C^\infty(P,K)^H$, we have 
$$(\Phi_\beta \circ \hat\alpha_1^*)(f)(t) 
= \beta(t).f(\hat\alpha_1(t))= f(\hat\alpha_1(t).\beta(t)^{-1}), $$
so that it suffices to see that the curves 
$\hat\alpha_0, \hat\alpha_1.\beta^{-1} \: [0,1] \to P$
with the same endpoints 
are homotopic with fixed endpoints, 
which is equivalent to the existence of a homotopy 
between $\hat\alpha_0.\beta$ and $\hat\alpha_1$. 

The homotopy $F$ can be lifted to a smooth map 
$\hat F \: [0,1] \times \R \to P$ such that the curves 
$\hat F_t$ are horizontal lifts starting in the base point $p_0$. 
Then $\hat F_i = \hat\alpha_i$ for $i =0,1$. 
If $\sharp$ denotes composition of paths and 
$\sim$ the homotopy relation, then 
the restriction of $\hat F$ to the boundary of $[0,1]^2$ 
shows that 
$\hat\alpha_1 
\sim \hat\alpha_0 \sharp \big(\hat\alpha_0(1).\beta\big) 
\sim \hat\alpha_0.\beta.$
\end{proof}

Putting all this information together, we now see 
with \eqref{eq:pull-rel}, \eqref{eq:pull-rel2} and \eqref{eq:pull-rel3} 
how $\cI_{\alpha_1} \circ \per_{\omega}$ and 
$\cI_{\alpha_0} \circ \per_{\omega}$ differ: 
\begin{align*}
\cI_{\alpha_1} \circ \per_{\omega} 
&= \per_{\omega_{\phi_1}} \circ \pi_2(\hat\alpha_1^*) 
= \per_{\L(\Phi_\beta)^*\omega_{\phi_0}} \circ \pi_2(\hat\alpha_1^*) 
- \per_{q_V \circ \eta_{\delta(\beta)}} \circ \pi_2(\hat\alpha_1^*) \\ 
&= \per_{\omega_{\phi_0}} \circ \pi_2(\Phi_\beta \circ \hat\alpha_1^*) 
- q_V \circ \per_{\eta_{\delta(\beta)}} \circ \pi_2(\hat\alpha_1^*) \\ 
&= \per_{\omega_{\phi_0}} \circ \pi_2(\hat\alpha_0^*) 
- q_V \circ \per_{\eta_{\delta(\beta)}} \circ \pi_2(\hat\alpha_1^*) \\ 
&= \cI_{\alpha_0} \circ \per_{\omega} 
- q_V \circ \per_{\eta_{\delta(\beta)}} \circ \pi_2(\hat\alpha_1^*).
\end{align*}

\begin{proposition} \label{prop:indep-crit} 
Let $\theta \in \Omega^1(P,\fh)$ be a principal 
connection form corresponding to the connection $\nabla$ and 
$R(\theta) = \dd\theta + \frac{1}{2} [\theta,\theta]\in \Omega^2(P,\fh)$ be 
its curvature. Then the homomorphisms 
$$ \cI_\alpha \circ \per_{\omega_\kappa^\nabla} \: 
\pi_2(\Gamma\cK) \to V_{\phi_V} $$
depend for each smooth loop $\alpha$ in $m_0$ only on the homotopy 
class $[\alpha] \in \pi_1(M,m_0)$ if and only if 
for each derivation $D \in \im(\L(\rho_\fk)\circ R(\theta))$,  
the periods of the cocycle $\eta_D(x,y) := \kappa(x,Dy)$ on $\fk$ 
are trivial. 
\end{proposition} 

\begin{proof} Suppose first that for 
each derivation $D \in \im(\L(\rho_\fk)\circ R(\theta))$,  
the periods of $\eta_D$ vanish. 
Let $F \: [0,1] \times \bS^1 \to M$ 
be a smooth homotopy of the loops $F_0$ and $F_1$ based in $m_0$. 
We  lift $F$ to a smooth map 
$\hat F \: [0,1]^2\to P$ such that the curves 
$\hat F_t = \hat F(t,\cdot)$ are horizontal, start in the base point $p_0$ 
and define the smooth curve $\beta \: [0,1] \to H$ by 
$\hat F_t(1) = \hat F_0(1).\beta(t)$. This implies that 
$$ \delta(\beta)_t = \theta(\frac{d}{dt} \hat F_t(1)) 
= (\hat F^*\theta)(\frac{\partial}{\partial t})(t,1). $$
Since the curves $\hat F_t$ are horizontal, 
$\hat F^*\theta(\frac{\partial}{\partial s})=0$, which leads to 
$$ \big(\hat F^*R(\theta)\big)\Big(
\frac{\partial}{\partial s},  \frac{\partial}{\partial t}\Big) 
= \dd(\hat F^*\theta)
\Big(\frac{\partial}{\partial s},  \frac{\partial}{\partial t}\Big) 
= \frac{\partial}{\partial s}(\hat F^*\theta)
\Big(\frac{\partial}{\partial t}\Big), $$
so that we arrive with 
$(\hat F^*\theta)(\frac{\partial}{\partial t})(t,0)=0$  
at 
\begin{equation}
  \label{eq:int-form-beta}
\delta(\beta)_t 
= (\hat F^*\theta)(\frac{\partial}{\partial t})(t,1) 
= \int_0^1 \big(\hat F^*R(\theta)\big)\Big(
\frac{\partial}{\partial s},  \frac{\partial}{\partial t}\Big)(t,s)\, ds. 
\end{equation}

We have to show that the periods of $\eta_{\delta(\beta)}$ vanish. 
As we have just seen, $\delta(\beta)_t = \int_0^1 D(s,t)\, ds$ 
holds for a smooth family $D(s,t)$ of derivations of $\fk$ for which 
the periods of $\eta_{D(t,s)}$ are trivial by assumption. 
Now the assertion follows 
from 
$$\eta_{\delta(\beta)} 
= \int_0^1 \eta_{\delta(\beta)_t}\, dt 
= \int_0^1 \int_0^1 \eta_{D(s,t)}\, ds\ dt, $$
so that 
$$ \per_{\eta_{\delta(\beta)}}([\sigma]) 
= \int_0^1 \int_0^1 \per_{\eta_{D(s,t)}}[\ev_t^K \circ \sigma]\, ds\ dt =0. $$

Now we prove the converse. 
To this end, we consider a family $(\alpha_t)_{0 \leq t \leq T}$ of 
smooth loops in $m_0$ with $\alpha_0 = m_0$ (the constant loop) 
for which the holonomy $\beta(t) = {\cal H}(\alpha_t)$ 
satisfies 
$$ \beta'(0) =0\quad \mbox{ and } \quad 
\beta''(0)= D := 2 R(\theta)_{p_0}(\tilde w, \tilde v) $$
for two horizontal vectors $\tilde v, \tilde w\in T_{p_0}(P)$. 
Since each loop $\alpha_t$ is contractible, the corresponding operator
$\phi_V$ on $V$ is $\id_V$ (cf.\ Appendix~\ref{app:D}). From
the assumption and
Remark~\ref{rem:4.3} we now obtain
that the cocycle 
$$ \eta(T) := \eta_{\delta(\beta)} = \int_0^T \eta_{\delta(\beta)_t}\, dt $$
has trivial periods for each $T > 0$. 
As a function of $T$, we have $\eta'(t) = \eta_{\delta(\beta)_t}$, 
$\eta'(0) = 0$ and $\eta''(0) = \eta_{\beta''(0)} = \eta_D$, so that 
all periods of $\eta_D$ must be trivial. 
\end{proof}

\begin{corollary} \label{cor:indep-crit} The homomorphism 
$\cI_\alpha \circ \per_{\omega}$ depends for each smooth loop 
$\alpha$ only on the homotopy class $[\alpha] \in \pi_1(M,m_0)$ 
if one of the following conditions is satisfied: 
\begin{itemize}
\item[\rm(a)] $\pi_2(K)$ is a torsion group. 
\item[\rm(b)] $\fk = \fh$ and $\rho_\fk = \ad$. 
\item[\rm(c)] The connection $\nabla$ on $P$ is flat. 
\end{itemize}
\end{corollary}

\begin{proof} (a) follows immediately from Proposition~\ref{prop:indep-crit} 
because $\tor \pi_2(K)$ lies in the kernel of each period homomorphism. 

\nin (b) For each inner derivation $D$ 
of $\fk$, the cocycle $\eta_D$ is a coboundary, so that 
its period map vanishes (\cite[Rem.~5.9]{Ne02a}). 
Therefore Proposition~\ref{prop:indep-crit} applies. 

\nin (c) Proposition~\ref{prop:indep-crit} 
applies because $R(\theta) =0$. 
\end{proof}

\begin{lemma} \label{lem:groupco} Let $\tilde m_0 \in \tilde M$ be a base 
point 
and realize $\Omega^1(M,\V)$ as the space 
$\Omega^1(\tilde M,V)^{\pi_1(M)}$ and identify $\pi_1(M)$ with the 
group of deck transformations of $\tilde M$, acting from the right. 
Then the map 
$$ \Psi \: H^1_{\rm dR}(M,\V) \to H^1(\pi_1(M),V), \quad 
\Psi([\theta]) = [\psi_\theta], \quad 
\psi_\theta(\gamma) := \int_{\tilde m_0}^{\tilde m_0.\gamma^{-1}} \theta$$
is a linear isomorphism. 
\end{lemma}

\begin{proof} The existence of the isomorphism $\Psi$ follows from 
\cite[p.~356]{CE56}. Here we briefly argue that it is given as above. 
Let 
$f_\theta \in C^\infty(\tilde M,V)$ be the unique function with 
$\dd f_\theta = \theta$ vanishing in $\tilde m_0$. 
From the observation that $\psi_\theta(\gamma) 
= \gamma .f_\theta - f_\theta$ is a constant function, it follows that 
$\psi_\theta$ is a $1$-cocycle and that $[\psi_\theta]$ only depends 
on the cohomology class $[\theta]$. 

If $\Psi([\theta]) = [\psi_\theta]$ vanishes, 
then there exists some $v \in V$ with 
$$ \psi_\theta(\gamma) = \gamma.f_\theta - f_\theta= \gamma.v - v. $$
Then $f_\theta- v$ is $\pi_1(M)$-invariant, so that 
$[\theta] = [\dd(f_\theta - v)] = 0$. Therefore $\Psi$ is injective. 

To see that $\Psi$ is also surjective, let $\chi \in Z^1(\pi_1(M),V)$. 
We consider the corresponding 
affine action of $\pi_1(M)$ on $V$, defined by 
$\gamma * v := \tilde\rho_V(\gamma).v - \chi(\gamma).$ 
Then the associated bundle $B := (\tilde M \times V)/\pi_1(M)$ is an affine 
bundle over $M$. Using smooth partitions of unity, we see that this 
bundle has a smooth section $s \: M \to B$. We write $s$ as 
$s(q_M(m)) = [(m, f(m))]$ for some smooth function 
$f \: \tilde M \to V$. Then 
$[(m.\gamma,f(m.\gamma))] = [(m,\gamma * f(m.\gamma))]$
implies 
$$ f(m) = \gamma * f(m.\gamma) = \tilde\rho_V(\gamma).f(m.\gamma) 
- \chi(\gamma) = (\gamma.f)(m)-\chi(\gamma), $$
so that $\chi(\gamma) = \gamma.f - f$  and thus 
$\chi = \psi_{\dd f}$. 
\end{proof}

\begin{remark} \label{rem:hausdorff} 
If $B^1(\pi_1(M),V)$ is a closed subspace of 
$Z^1(\pi_1(M),V)$ with respect to the topology of pointwise convergence, 
then it follows from the proof of the Lemma~\ref{lem:groupco} 
that $\dd\big(C^\infty(\tilde M,V)^{\pi_1(M)}\big)$ 
is a closed subspace 
of $\Omega^1(\tilde M,V)^{\pi_1(M)}$, 
but in general this is not clear. 

If $\dim V <~\infty$, then $B^1(\pi_1(M),V)$ is finite-dimensional, hence 
closed, so that the quotient space $H^1(\pi_1(M),V) \cong 
H^1_{\rm dR}(M,\V)$ is Hausdorff. 
\end{remark}

\begin{remark} There are compact $3$-manifolds $M$ 
whose fundamental group contains normal subgroups which are not finitely 
generated. In fact, by van Kampen's Theorem, 
the fundamental group of the connected sum $M$ of four 
copies of $\bS^1 \times \bS^2$ is the free group on $4$ generators. 
In VA.23(iv) of \cite{dlH00} one finds an 
example of a group with four generators which is not finitely 
presented. This implies that the normal subgroup $R$ of $\pi_1(M)$, 
the free group 
on four generators, generated by these relations is not finitely 
generated. Therefore $\tilde M/R$ is a $3$-manifold whose fundamental 
group $R$ is not finitely generated, although $\tilde M/R$ covers 
the compact manifold $M$. 
\end{remark}

{\sloppy
\begin{proposition}\label{prop:reduction} 
For each \mbox{$\alpha \in C^{\infty}(\bS^{1},M)$}, 
the homomorphism \mbox{$\cI_\alpha \circ \per_{\omega}$} 
depends only on the homotopy class \mbox{$[\alpha] \in \pi_1(M,m_0)$} 
if and only if \mbox{$\Pi_{\omega} \subeq H^1_{\rm dR}(M,\V)$}. 
\end{proposition}}

\begin{proof} We realize 
$\Omega^1(M,\V)$ as the space $\Omega^1(\hat M,V)^\D$ of 
$\D$-invariant $V$-valued $1$-forms on $\hat M$. 
Let $\hat\alpha_M$ denote the image of the horizontally lifted curve
$\hat\alpha \: [0,1] \to P$ in $\hat M = P/\ker \rho_V$ 
and observe that for two homotopic loops $\alpha_0$ and $\alpha_1$ in $m_0$, 
the curves $\hat\alpha_{1,M}$ and $\hat\alpha_{2,M}$ are 
homotopic with fixed endpoints. 

If, conversely, $\beta_0,\beta_1 \: [0,1] \to \hat M$ 
are two smooth curves starting in $\hat m_0$ which are homotopic 
with fixed endpoints, such that the curves 
$\alpha_i := \hat q_M \circ \beta_i$ are closed, then 
$\beta_i = \hat\alpha_{i,M}$ for $i =1,2$. 

Since $\hat M \cong P/\ker \rho_V$ is connected, the 
homomorphisms $\cI_\alpha \circ \per_{\omega}$ 
depend only on the homotopy class of $\alpha$ 
if and only if each element 
$[\theta] \in \Pi_{\omega} \subeq \oline\Omega^1(M,\V) 
\subeq \Omega^1(\hat M,V)^\D$ has the property that the 
integral $\cI_\alpha([\theta]) = \int_{\hat\alpha_M} \theta$ 
only depends on the homotopy class of $\hat\alpha_M$ 
with fixed endpoints. This is equivalent to the $1$-form $\theta$ 
being closed, i.e., $[\theta] \in H^1_{\rm dR}(M,\V)$. 
\end{proof}

\begin{remark} \label{rem:Dfin} 
If the group $\D$ is finite, then the fixed point functor
$H^0(\D,\cdot)$ is exact on rational $\D$-modules, so that
\begin{align*}
 H^1_{\rm dR}(M,\V) 
&= Z^1_{\rm dR}(\hat M,V)^\D/\dd\big(C^\infty(\hat M,V)^\D\big)\\
&= Z^1_{\rm dR}(\hat M,V)^\D/B^1_{\rm dR}(\hat M,V)^\D 
= H^1_{\rm dR}(\hat M,V)^\D. 
\end{align*}
Since $\D \cong \pi_1(M)/\pi_1(\hat M)$ is finite 
and $V$ is divisible, the surjective map 
$$ H^1_{\rm dR}(M,V) \cong  \Hom(\pi_1(M),V) \to 
 H^1_{\rm dR}(\hat M,V) \cong  \Hom(\pi_1(\hat M),V) $$ 
is a linear isomorphism, and we thus obtain 
\begin{equation}
  \label{eq:derham-fin}
H^1_{\rm dR}(M,\V) 
\cong H^1_{\rm dR}(\hat M,V)^\D 
\cong H^1_{\rm dR}(M,V)^\D \cong H^1_{\rm dR}(M,V^\D). 
\end{equation}
\end{remark}

\begin{remark} For $M = \bS^1 \cong \R/\Z$ and
the $m$-fold covering $\hat M := \R/m\Z \cong \bS^1$, we have
$\D \cong \Z/m$ and 
$H^1_{\rm dR}(M,\V) \cong H^1_{\rm dR}(M,V)^\D \cong V^\D.$
\end{remark}

\begin{theorem}\label{thm:reduction} {\rm(Reduction Theorem)} 
Assume that $\D$ is finite and that $\Pi_{\omega} \subeq H^1_{\rm dR}(M,\V)$. 
Then $\Pi_{\omega}$ is discrete if this is 
the case for each $\Pi_{\omega_{\alpha }}$,
with $\omega_{\alpha}:=\omega_{\varphi_{K}^{\alpha}}$ for $\alpha \in
C^\infty(\bS^1,\hat M)$.
\end{theorem}

\begin{proof} Since $\D \cong \pi_1(M)/\pi_1(\hat M)$ is finite, 
there exists a number $N \in \N$ such that the image of the 
homomorphism $H_1(\hat M) \to H_1(M)$ contains $N\cdot H_1(M)$. 
Suppose that $H_{1}(M)$ is finitely generated of rank $r$. 
Then the Universal Coefficient Theorem, combined with de 
Rham's Theorem, yields  
$$ H^1_{\rm dR}(M,V) \cong \Hom(H_1(M),V) \cong V^r. $$
As we have seen above, there exist smooth loops 
$\oline\alpha_i \in C^{\infty}_{*}(\bS^{1},\hat M)$, $i =1,\ldots, r$, 
whose images in $H_1(M)$ form a $\Q$-basis of  
$H_{1}(M)\otimes \Q$. We then obtain the concrete linear isomorphism 
\[
\Phi = (\cI_{\alpha_i})_{i = 1,\ldots, r}\from H^{1}_{\rm dR}(M,V) 
\to V^r, \quad [\omega]\mapsto 
\Big(\int_{\ol{\alpha}_{i}}\omega\Big)_{i=1,\ldots, r}. 
\]
By \eqref{eq:pull-rel2},
$\Phi(\Pi_\omega) \subeq \prod_{i =1}^r \Pi_{\omega_{\alpha_i}},$
where the right hand side is a discrete subgroup of $V^r$. 
Therefore $\Pi_\omega$ is discrete. 
\end{proof}

\begin{remark}
The assumption on $\D$ to be finite in the previous theorem was
needed to ensure that the map 
$$\Phi \: H^1_{\rm dR}(M,\V) \to H^1_{\rm dR}(\hat M,V)^\D$$ 
is an isomorphism. The argument also works if 
$\Phi$ is injective and $H_1(\hat M)$ is finitely 
generated. The kernel of 
$\Phi \: H^1(\pi_1(M),V) \to H_1(\pi_1(\hat M),V)$ is the image of the 
natural map $H^1(\D,V) \to H^1(\pi_1(M),V)$, hence vanishes whenever
$H^1(\D,V)=\{0\}$. 

For a finite-dimensional orthogonal representation of $\D$ on $V$,
this is the case if $\D$ has Kazhdan's property (T) (\cite[Prop.~7 and
Prop.~31]{Pa07}).
\end{remark}

Combining the Reduction Theorem with Corollary~\ref{cor:indep-crit}, 
leads to: 

\begin{corollary} If $\D$ is finite, $K = H$ and $\fK = \Ad(P)$, 
then $\Pi_{\omega}$ is discrete if this is 
the case for each $\Pi_{\omega_{\alpha}}$, 
$\alpha \in C^\infty(\bS^1,\hat M)$. 
\end{corollary}

\begin{theorem} \label{thm:redux-a} If $\pi_2(K)$ vanishes, then 
the following are equivalent:

\begin{description}
\item[\rm(1)] $\Pi_{\omega}$ is discrete for each compact manifold $M$ 
and each connection $\nabla$, provided the group $\D \cong \rho_V(H)$ 
is finite. 
\item[\rm(2)] $\Pi_{\omega}$ is discrete for the trivial bundle over 
$\bS^1$ and the canonical connection. 
\item[\rm(3)] The period group $\Pi_{C(\kappa)}$ of the $3$-cocycle 
$C(\kappa)$  of $\fk$ is discrete. 
\end{description}
\end{theorem}

\begin{proof} The equivalence of (2) and (3) follows from 
Lemma~\ref{lem:per-loop}, so that it remains to derive (1) from (2). 
If $\pi_2(K)$ vanishes, then 
Lemma~\ref{lem:per-loop} further implies that 
the period group of any cocycle 
$\omega_\alpha$, $\alpha \in C^\infty(\bS^1,\hat M)$, 
is discrete if and only if this is the case for $\Pi_{C(\kappa)}$. 
Now the Reduction Theorem~\ref{thm:reduction} applies. 
\end{proof}

With Corollary~\ref{cor:disc-univ} we also get: 
\begin{corollary}\label{cor:disc-univ-general}
If $\fk$ is finite-dimensional, $V = V(\fk)$, $\kappa = \kappa_u$ 
is universal, and $\D$ is finite, then the period group of the cocycle 
$\omega_\kappa^\nabla$ is discrete for any connection $\nabla$.
\end{corollary}

\msk

\begin{remark} If $K$ is finite-dimensional and $1$-connected, then 
$H := \Aut(\fk) \break \cong \Aut(K)$ has finitely many connected components 
because $\Aut(\fk)$ is a real algebraic group (\cite{OV90}). 

If $V^0(\fk) 
:= V(\fk)/\der(\fk).V(\fk)= V(\fk)/\der(\fk).V_0(\fk)$ 
denotes the quotient space, then 
the corresponding form $\kappa_u^0 \: \fk \times \fk \to V^0(\fk)$ 
is the universal $\der(\fk)$-invariant symmetric bilinear form. 
Then $\kappa_u^0$ 
is invariant under $\Aut(\fk)_0$ and $\pi_0(\Aut(\fk))$ is finite.
Since the period groups of $C(\kappa_u)$ and $C(\kappa_u^0)$ coincide 
(Theorem~\ref{thm:disc1}), Theorem~\ref{thm:redux-a} 
implies that the period group 
$\Pi_{\omega_{\kappa_u^0}}$ is discrete. 
\end{remark}

\begin{example} \label{ex:counter-ex-h1} 
Now we show that $\Pi_\omega$ is not always contained in 
$H^1_{\rm dR}(M,\V)$. 

We consider a trivial bundle $\cK = M \times K$ and $H = \R$, so that 
$\fh = \R$ acts on $\fk$ by a derivation $D \in (\der\fk)_\kappa$ 
and  the bundle $\V$ is trivial. 
We then write any covariant exterior derivative as 
$$ \dd^\nabla f = \dd f + \beta\cdot Df, \quad f \in \Omega^1(M,\R) $$
for some $\beta \in \Omega^1(M,\R)$, 
and, accordingly, $\omega = \omega_0 + \eta_\beta 
= \omega_0 + \beta \cdot \eta_D$. 
Since all periods of $\omega_0$ are contained in $H^1_{\rm dR}(M,V)$ 
(Corollary~\ref{cor:indep-crit}(d)), $\Pi_\omega$ is contained 
in $H^1_{\rm dR}(M,\V)$ if and only if this holds for 
the period group of $\eta_\beta$. 

On $\bS^1$, each $1$-form is closed, so that we consider 
$M = \T^2$. Then the range of 
$$ \per_{\eta_\beta} = \beta \cdot \per_{\eta_D} 
\: \pi_2(K) \to \Omega^1(M,V) $$
does not lie in the space of closed forms if 
$\beta$ is not closed and $\per_{\eta_D}$ is non-trivial, 
which is the case for $K = C^\infty(\bS^1,\g)$, 
$\g$ simple compact and $Df = f'$ 
(cf.\ Example~\ref{ex:torus}).
\end{example}

\subsection{Integrating actions}

In this section we show that for any principal $K$-bundle $P$ 
($K$ locally exponential), 
the action of the Lie group $\Aut(P)$ of bundle automorphism 
on the spaces $\oline\Omega^1(M,\V)$ and the Lie algebra 
$\gau(P)$ combines to a smooth automorphic action on the central 
Lie algebra extension $\hat\gau(P)$, defined by the cocycle 
$\omega = \omega_\kappa^\nabla$. 
Moreover, we show under which conditions this construction carries
over to arbitrary Lie group bundles which are not necessarily 
gauge bundles.

Let $\theta \in \Omega^1(P,\fk)$ be a principal connection $1$-form 
corresponding to $\nabla$. 
Realizing $\gau(P)$ in $C^\infty(P,\fk)$, 
we have $\nabla f = \dd f + [\theta,f]$, so that 
$$ \omega(f_1, f_2) 
= [\kappa(f_1, \nabla f_2)] 
= [\kappa(f _1, \dd f_2) + \kappa(\theta,[ f_2,  f_1])]. $$
The Lie group $\Aut(P)$ acts smoothly on the affine space 
${\cal A}(P) \subeq \Omega^1(P,\fk)$ of principal connection 
$1$-forms by $\phi.\theta := (\phi^{-1})^*\theta$ and on 
$\gau(P)$ by $\phi. f =  f \circ \phi^{-1}$ 
(cf.\ \cite[Prop.~6.4]{Gl06}). 
We then have 
$$ \phi.\dd^\nabla f 
= \phi.(\dd  f + [\theta, f])
= \dd(\phi. f) + [\phi.\theta,\phi. f])
= \dd^\nabla(\phi. f) + [\phi.\theta-\theta,\phi. f]. $$
This leads to 
$$ (\phi.\omega)( f_1,  f_2) 
=  \phi.\omega(\phi^{-1}. f_1, \phi^{-1}. f_2) 
=  \omega( f_1,  f_2) + [\kappa(\phi.\theta-\theta,[ f_2,  f_1])].$$
Note that 
$$\zeta \: \Aut(P) \to \Omega^1(M,\Ad(P)),\quad \phi \mapsto 
\phi.\theta -\theta $$
is a smooth $1$-cocycle, so that 
$$ \Psi \: \Aut(P) \to \Hom(\gau(P),\oline\Omega^1(M,\bV)), \quad 
\Psi(\phi)( f) := [\kappa(\phi.\theta -\theta, f)] $$
is a $1$-cocycle with 
$\dd_{\gau(P)}(\Psi(\phi)) = \phi.\omega- \omega$, defining a smooth map 
$$\Aut(P) \times \gau(P) \to \oline\Omega^1(M,\V).$$

\begin{theorem}\label{thm:aut(P)-action} 
The group $\Aut(P)$ acts smoothly by automorphisms on the centrally 
extended Lie algebra $\hat\gau(P)$ by 
\begin{equation}
  \label{eqn:aut(P)action} 
\phi.([\alpha], f) := ([\phi.\alpha] 
+ [\kappa(\phi.\theta-\theta,\phi. f)], \phi. f).
\end{equation}

If, in addition, the period group $\Pi_\omega$ 
is discrete and $Z\hookrightarrow \wh{G}\twoheadrightarrow G$ is a 
central extension with $1$-connected $\wh{G}$, $G=\Gau(P)_0$ and 
Lie algebra $\hat\g$, then 
this action integrates to a smooth action of $\Aut(P)$ on $\hat G$.
\end{theorem}

\begin{proof} First, \cite[Lemma~V.1]{MN03} implies that we 
obtain automorphisms of $\hat\gau(P)$, and the smoothness of the 
action follows from the smoothness of $\zeta$ and the smoothness 
of the actions of $\Aut(P)$ on $\gau(P)$ and $\oline\Omega^1(M,V)$. 

Assume that $\Pi_\omega$ is discrete. 
Since the action of $\Aut(P)$ on $\g := \gau(P)$ and 
$\z := \oline\Omega^1(M,\V)$ preserves the cohomology class of $\omega$ (cf.\ Example \ref{ex:1.4}), 
the period homomorphism 
$\per_\omega \: \pi_2(G) \to \z$
is $\Aut(P)$-equivariant, which implies in particular that its image  
in $\z$ is invariant under the action of $\Aut(P)$. We therefore 
obtain a smooth action of $\Aut(P)$ on $Z_0 := \z/\Pi_\omega$. 
Now the group $\hat G$ is a central extension of the 
universal covering group $\tilde G$ of $G$ by $Z_0$, and 
$\pi_{0}(Z)\cong \pi_1(G)$ 
(cf.\ \cite[Rem.~7.14]{Ne02a}). 
Finally, we lift the $\Aut (P)$ action on $G$ to a smooth action on 
$\tilde G$ and apply the Lifting Theorem \ref{thm:lifting-theorem}.
\end{proof}

If $\phi_f \in \Gau(P)$ is a gauge transformation corresponding 
to the smooth function $f \: P \to K$, then 
$\phi_f^*\theta = \delta(f) + \Ad(f)^{-1}\theta$ implies 
$$ \phi_f.\theta = \delta(f^{-1}) + \Ad(f)\theta 
\quad \mbox{ and } \quad  
\zeta(\phi_f) = \delta(f^{-1}) + \Ad(f)\theta -\theta. $$

\begin{corollary}\label{cor:int-finite-dim}
The adjoint action of $\gau(P)$ on $\hat\gau(P)$ integrates to a smooth 
action of $\Gau(P)$ on $\hat\gau(P)$. 
\end{corollary}

\begin{theorem} \label{thm:redux} If $\pi_{0}(K)$ is finite and 
$\pi_2(K)$ vanishes, then the following are
equivalent:
\begin{description}
\item[\rm(1)] $\omega_{\kappa}^\nabla$ integrates for each principal 
$K$-bundle $P$ over a compact manifold $M$ to a Lie group extension of 
$\Gau(P)_0$. 
\item[\rm(2)] $\omega_{\kappa}$ integrates for the trivial 
$K$-bundle $P = \bS^1 \times K$ over $M = \bS^1$ 
to a Lie group extension of $C^\infty(\bS^1,K)_0$. 
\item[\rm(3)] The image of $\per_{\kappa} \: \pi_3(K) \to V$ is discrete. 
\end{description}
\end{theorem}

\begin{proof}
Since the existence of a Lie group extension of $G := \Gau(P)_0$ integrating
$\omega_{\kappa}^{\nabla}$ is equivalent to the discreteness of
$\Pi_{\omega_{\kappa}^{\nabla }}$ and the integrability of the adjoint
action to an action on $\hat\gau(P)$, this follows from 
Corollary~\ref{cor:int-finite-dim} and Theorem~\ref{thm:redux-a}.
\end{proof}

\begin{theorem} 
Let $P$ be a finite-dimensional connected principal bundle with structure group
$K$ over the compact manifold $M$. If $V = V(\fk)$, $\kappa =\kappa_u$ 
is universal and $\D = \oline\rho_{V}(\pi_0(K)) \subeq \GL(V(\fk))$ 
is finite, then the central extension $\hat\gau(P)$ of $\gau(P)$ 
defined by $\omega_{\kappa}^{\nabla}$ 
integrates for any connection $\nabla$ to a central extension of 
the identity component $\Gau(P)_{0}$ of the gauge group.
\end{theorem}

\begin{proof}
With Corollary~\ref{cor:disc-univ-general}, this follows as in the proof of
Theorem~\ref{thm:redux}.
\end{proof}

For general Lie algebra bundles $\fK$, the action of 
$\Gamma\fK$ on $\hat{\Gamma\fK}$ is given by 
$$ h.([\alpha], f) 
= (\omega(h,f), [h,f]) 
= ([\kappa(-\dd^\nabla h,f)], [h,f]). $$
The fact that $\nabla$ is a Lie connection means that 
$$ \dd^\nabla \: \Gamma\fK \to \Omega^1(M,\fK) $$
is a $1$-cocycle for the action of the Lie algebra $\Gamma\fK$ on 
$\Omega^1(M,\fK)$ by $f.\alpha := [f,\alpha]$ (pointwise bracket). 
To integrate the action of $\Gamma\fK$ on $\hat{\Gamma\fK}$ to a 
group action, we therefore have to integrate $\dd^\nabla$ 
to a Lie group cocycle 
$(\Gamma\fK)_0 \to \Omega^1(M,\fK).$

This can be achieved as follows. We assume that 
$K$ is $1$-connected. First we observe that 
$\der(\fk) = Z^1(\fk,\fk)$, where $\fk$ acts on itself by the 
adjoint action. In this sense, each derivation $D \in \der(\fk)$ 
is a $1$-cocycle, hence defines an equivariant closed $1$-form 
$D^{\rm eq} \in \Omega^1(K,\fk)$ which is exact since $K$ 
is $1$-connected. Let $\chi^D \: K \to \fk$ be the unique 
smooth function with $\dd \chi^D = D^{\rm eq}$ and $\chi^D(\1) = 0$. 
Then $\chi^D$ is a smooth $1$-cocycle (cf.~\cite{GN09}), 
and the smoothness of the 
action $\rho_\fk$ of $\h$ on $\fk$ implies that the function 
$$ \chi \: \h \times K \to \fk, \quad (x,k) \mapsto 
\chi^{\rho_\fk(x)}(k) $$
is smooth. 

If $\theta \in \Omega^1(P,\fh)$ is a principal connection 
$1$-form, we now define for $f \in \Gamma\cK$ 
a $1$-form $\chi^\theta(f)$ in $\Omega^1(P,\fk)^H \cong \Omega^1(M,\fK)$ by 
$$ \chi^\theta(f)v := \chi^{\theta(v)}(f(p)) \quad \mbox{ for } \quad 
v \in T_p(P). $$
Then 
$$\delta^\nabla(f) := \delta(f) + \chi^{\theta}(f^{-1}) $$
is a covariant left logarithmic derivative on $\Gamma\cK$ and 
$$ \Gamma\cK \to \Omega^1(M,\fK), \quad 
f \mapsto \delta^\nabla(f^{-1}) $$
is a $1$-cocycle integrating $-\dd^\nabla$.
We now calculate for $\phi \in \Gamma\cK$:  
\begin{align*}
\phi.\dd^\nabla f 
&= \phi.(\dd  f + \theta.f)
= \dd(\phi.f) + (\phi.\theta).(\phi.f) \\
&= \dd^\nabla(\phi.f) + \big(\phi.\theta-\theta).(\phi.f)
= \dd^\nabla(\phi.f) + \chi^\theta(\phi).(\phi.f).
\end{align*}
This easily leads to 
$$ (\phi.\omega)( f_1,  f_2) 
=  \phi.\omega(\phi^{-1}. f_1, \phi^{-1}. f_2) 
=  \omega( f_1,  f_2) + [\kappa(\chi^\theta(\phi),[ f_2,  f_1])],$$
and 
$\chi^\theta \: \Gamma\cK  \to \Omega^1(M,\fK)$ 
is a smooth $1$-cocycle, so that 
$$ \Psi \: \Gamma\cK \to \Hom(\Gamma\fK,\oline\Omega^1(M,\bV)), \quad 
\Psi(\phi)( f) := [\kappa(\chi^\theta(\phi), f)] $$
is a $1$-cocycle with 
$\dd_{\Gamma\fK}(\Psi(\phi)) = \phi.\omega- \omega$, defining a smooth map 
$$\Gamma\cK\times \Gamma\fK \to \oline\Omega^1(M,\V).$$ 

\begin{theorem}\label{thm:fK-action} 
If $K$ is $1$-connected, then the group $\Gamma\cK$ 
acts smoothly by automorphisms on the centrally 
extended Lie algebra $\hat{\Gamma\fK}$ by 
\begin{equation}
  \label{eqn:fK-action} 
\phi.([\alpha], f) := ([\phi.\alpha] 
+ [\kappa(\chi^\theta(\phi),\phi. f)], \phi. f).
\end{equation}
\end{theorem}

\begin{theorem} \label{thm:redux-fK} If 
$K$ is $2$-connected, then the following are equivalent:
\begin{description}
\item[\rm(1)] If $\rho_V(H)$ is finite, 
then the Lie algebra defined by $\omega_{\kappa}^\nabla$ 
integrates for each $K$-bundle $\cK$ over a compact manifold $M$ 
and each Lie connection $\nabla$ on $\fK$ 
to a Lie group extension of $(\Gamma\cK)_0$. 
\item[\rm(2)] The extension defined by 
$\omega_{\kappa}$ on $C^\infty(\bS^1,\fk)$ 
integrates to a Lie group extension of $C^\infty(\bS^1,K)_0$. 
\item[\rm(3)] The image of $\per_{\kappa} \: \pi_3(K) \to V$ is discrete. 
\end{description}
\end{theorem}

\begin{proof}
Since the existence of a Lie group extension of $(\Gamma\cK)_0$ 
integrating 
$\omega_{\kappa}^{\nabla}$ is equivalent to the discreteness of
$\Pi_{\omega_\kappa^\nabla}$ and the integrability of the adjoint
action of $\Gamma\fK$ on $\hat{\Gamma\fK}$ to an action of $(\Gamma\cK)_0$, 
this follows from Theorem~\ref{thm:fK-action} and Theorem~\ref{thm:redux-a}.
\end{proof}
\appendix

\section{Appendix: The Lie group structure on $\Gamma{\cal K}$} \label{app:A} 

In this section we explain how to obtain a locally exponential Lie group 
structure on the group $\Gamma{\cal K}$ of smooth sections of the 
(locally trivial) Lie group bundle ${\cal K}$ over the compact manifold 
$M$ whose fiber is a locally exponential Lie group $K$ with Lie algebra $\fk$. 

We further assume that the Lie group bundle $\cK$ is associated to a 
principal $H$-bundle $P$ via a smooth action defined by
$\rho_K \: H\to \Aut (K)$. We write $\rho_\fk(h) := \L(\rho_K(h))$ 
for the corresponding smooth action of $H$ on $\fk$ 
and $\fK := \L({\cal K})$ for its Lie algebra 
bundle with fiber $\fk :=\L(K)$. 
We endow the space $\Gamma\fK$ of smooth sections of $\fK$
with the smooth compact open  topology.  This turns $\Gamma\fK$ into a
locally convex Lie algebra because over each open subset $U \subeq M$
for which ${\cal K}_U$ is trivial, we have $\Gamma(\fK_U) \cong
C^\infty(U,\fk)$, and the Lie bracket on the locally convex space
$C^\infty(U,\fk)$ is continuous since $U$ is finite-dimensional.  
Likewise, the smooth compact open topology turns the
group $\Gamma{\cal K}$ of smooth section of ${\cal K}$ into a
topological group. Indeed, restriction 
defines a group homomorphism 
$\Gamma\cK \to \Gamma\cK_U \cong C^\infty(U,K)$, and the topology 
on $\Gamma\cK$ is defined by the embedding 
$\Gamma\cK \into \prod_U C^\infty(U,K)$, where $U$ runs through 
an open cover of $M$ consisting of trivializing open subsets 
(cf.\ \cite[Def.~II.2.7]{Ne06}). 

Since the exponential function $\exp_K \: \fk \to K$ is natural, we have 
$$ \exp_K \circ \L(\phi) = \phi \circ \exp_K $$
for every automorphism $\phi \in \Aut(K)$, and we obtain a fiberwise 
defined exponential map 
$$ \exp_{\cal K} \: \fK \to {\cal K}. $$
Composing smooth sections with this exponential map, we obtain a map 
$$ \exp_{\Gamma{\cal K}} \: \Gamma\fK \to \Gamma{\cal K}. $$

\begin{theorem} \label{thm:liestruc} 
The topological group $\Gamma{\cal K}$ carries a locally 
exponential Lie group structure with 
$\L(\Gamma{\cal K}) \cong \Gamma\fK$. 
 Moreover, this topology coincides with the smooth compact-open topology. 
\end{theorem}

\begin{proof} The proof of \cite[Thm.~1.11]{Wo07a} carries over from
the case of the conjugation of $K$ on itself to an arbitrary action of 
some group $H$ on $K$.
\end{proof}

\section{Appendix: The universal invariant bilinear form in finite dimensions} 
\label{app:B} 

Throughout this section, $K$ denotes a finite-dimensional Lie group 
and $\fk = \L(K)$ its Lie algebra. We further choose a Levi decomposition 
$\fk = \fr \rtimes \fs$ and write 
$\fs = \fs_1^{m_1} \oplus \ldots \oplus \fs_r^{m_r},$
for the decomposition of $\fs$ into simple ideals $\s_i$,  
where $\s_i$ is supposed to be non-isomorphic to $\s_j$ for $j \not=i$. 

\subsection{The action of $\Aut(\fk)$ on  $V(\fk)$} 

\begin{definition}
We put $V(\fk) := S^2(\fk)/\fk.S^2(\fk)$, 
where the action of $\fk$ on $S^2(\fk)$ is the natural action 
inherited by the one on the 
tensor product $\fk \otimes \fk$ by 
$x.(y \otimes z) = [x,y] \otimes z + y \otimes [x,z]$. 
There exists a natural invariant symmetric bilinear form 
$$\kappa_u \: \fk \times \fk \to V(\fk), \quad (x,y) \mapsto [x \vee y]$$ 
such that for each 
invariant symmetric bilinear form $\beta \: \fk \times \fk \to W$ there
exists a unique linear map $\phi \: V(\fk) \to W$ with $\phi \circ
\kappa_u = \beta$. We call $\kappa_u$ the 
{\it universal invariant symmetric bilinear form on $\fk$}. 
\end{definition}

\begin{remark} \label{rem:redu} (cf.\ \cite{MN03}) 
(a) The 
assignment $\g \to V(\g)$ is a covariant functor from finite-dimensional Lie
algebras to vector spaces. 

\par\nin(b) If $\g = \fa \oplus \fb$ and $\fa$ is perfect, then 
$V(\g) \cong V(\fa) \oplus V(\fb)$ because for every symmetric invariant
bilinear map $\kappa \: \g \times \g \to V$, we have for $x,y \in \fa$,
$z \in \fb$ the relation $\kappa([x,y],z) 
= \kappa(x,[y,z]) = \kappa(x,0) = 0$. 

\par\nin(c) If $\fh \trile \g$ is an ideal and 
the quotient morphism $q \: \g \to \fq := \g/\fh$ splits, then 
$\g \cong \fh \rtimes \fq$, and the natural map $V(\fq) \to V(\fg)$ is an
embedding. In fact, let $\eta \: \fq \to \g$ be the inclusion map. 
Then $q \circ \eta = \id_\fq$ and 
this leads to $V(q) \circ V(\eta) = \id_{V(\fq)}$, showing that
$V(\eta)$ is injective. 

\par\nin(d) If $\fs$ is reductive with the simple ideals $\fs_1,\ldots,
\fs_n$, then (b) implies that 
$$V(\fs) \cong V(\z(\fs)) \oplus \bigoplus_{j = 1}^n V(\fs_j).$$ 

\par\nin(e) If $\fk = \fr \rtimes \fs$ is a Levi decomposition, then (c)
shows that the natural map $V(\fs) \to V(\fk)$ is an embedding. 
\end{remark}

\begin{remark} \label{rem:aut-act} 
As a consequence of our construction, 
the group $\Aut(\fk)$ and its Lie algebra $\der(\fk)$ act 
naturally on $V(\fk)$. The Lie algebra $\fk$ itself, resp., the subalgebra 
$\ad \fk$ of inner derivations acts trivially. 

If all derivations are inner, as is the case if $\fk$ is semisimple, 
it follows that the identity component $\Aut(\fk)_0$ acts trivially 
on $V(\fk)$. 
\end{remark}

\begin{remark} \label{rem:simple} For 
a simple finite-dimensional real Lie algebra $\fs$, its {\it centroid}
$$ \Cent(\fs) := \{ \phi \in \End(\fs) \: (\forall x \in \fs)\ 
[\phi,\ad x] =0\} $$
is a field, hence isomorphic to $\R$ or $\C$ (\cite[Theorem~X.1]{Ja79}). 
If $\Cent(\fs) \cong \C$, then $\fs$ actually carries the structure of a 
complex simple Lie algebra and if $\Cent(\fs) \cong \R$, then 
its complexification $\fs_\C$ is simple.  In the latter case we call $\fs$ \textit{central simple}. 

If $\beta(x,y) := \tr(\ad x\ad y)$ is the Cartan--Killing form of $\fs$, 
then the map 
$$ \eta \: \Cent(\fs) \to \Sym^2(\fs,\R)^\fs \cong V(\fs)^*, \quad 
\eta(\phi)(x,y) := \beta(\phi(x),y) $$
is easily seen to be a linear isomorphism. This implies that 
$V(\fs)\cong \C$ if $\fs$ is complex and $V(\fs) \cong \R$ otherwise. 
In the latter case the Cartan--Killing form $\beta$ is already universal, 
and in the former case, we have the additional form 
$\beta'(x,y) = \beta(ix,y)$. Hence the Cartan--Killing form 
$$\beta_\C \: \fs \times \fs \to \C, \quad 
\beta_\C(x,y) 
= \frac{1}{2}\big(\beta(x,y) - i \beta(ix,y)\big) 
= \frac{1}{2}\big(\beta(x,y) - i \beta(x,iy)\big) $$ 
of the complex simple Lie algebra 
$\fs$ is the universal invariant symmetric bilinear form for the real simple 
Lie algebra $\fs$. 
\begin{footnote} {If $C$ is a complex linear endomorphism of a complex 
vector space, then the traces of $C$ with respect to $\R$ and $\C$ are related by 
$\tr_\C C = \frac{1}{2}(\tr_\R C - i \tr_\R(iC)).$}
\end{footnote}
\begin{footnote}{\cite{MN03}, Remark~II.2(4) 
uses the invalid assumption that $V(\fs)$ is one-dimensional 
for any real simple Lie algebra 
$\fs$. This has no serious consequence for the validity of the main results in 
that paper. The corresponding gap in the proof of Theorem~II.9 loc.\ cit.\ 
is fixed by Proposition~\ref{prop:simpcase} and Theorem~\ref{thm:disc1} below. 
Moreover, the assertion of Lemma~II.11 loc.\ cit.\ should read 
$V(\fk \otimes A) \cong V(\fk) \otimes A$ for $\fk$ simple finite-dimensional.}
\end{footnote}
\end{remark}

\begin{example}
If $\fk = \fgl_n(\R)$, then Remark~\ref{rem:redu}(d) implies that 
$$V(\fgl_n(\R)) 
\cong V(\fsl_n(\R)) \oplus V(\R) \cong \R^2 $$
because $\fsl_n(\R)$ is central simple. 
\end{example}

\begin{theorem} \label{thm:decomp} 
Let $\fk$ be a finite-dimensional real Lie algebra 
with Levi decomposition $\fk = \fr \rtimes \fs$ and 
$\fs = \bigoplus_{i = 1}^r \fs_i^{m_i}$ the decomposition into simple ideals. 
With $V_0 := \kappa_u(\fr, \fk)$ and 
$V_i := \kappa_u(\fs_i^{m_i},\fs_i^{m_i}) \cong V(\fs_i)^{m_i}$, 
we obtain a direct sum decomposition 
\begin{equation}
  \label{eq:deco}
V(\fk) = V_0 \oplus V_1 \oplus \ldots \oplus V_r 
\cong  V_0 \oplus V(\fs_1)^{m_1} \oplus \ldots \oplus V(\fs_r)^{m_r}
\end{equation}
which is invariant under the group $\Aut(\fk)$. 
\end{theorem} 

\begin{proof} Let $q \: \fk \to \fk/\fr \cong \fs$ denote the 
quotient map. Then Remark~\ref{rem:redu}(c) implies that 
$V(\fs)$ can be identified with a complement of the kernel of 
$V(q)$. Clearly, $\ker V(q) \supeq \kappa_u(\fr,\fk)$, 
and since 
$$V(\fk) 
= \kappa_u(\fk,\fk) 
= \kappa_u(\fr,\fk) + \kappa_u(\fs,\fs) 
= V_0 + V(\fs), $$
we see that $\ker V(q) = V_0$ and that the sum of $V(\fs)$ and 
$V_0$ is direct. 
The decomposition of $V(\fs)$ follows from Remark~\ref{rem:redu}(d). 

Now we show that the decomposition \eqref{eq:deco} is invariant 
under $\Aut(\fk)$. 
Let $\Inn(\fk) \subeq \Aut(\fk)_0$ denote the normal subgroup of 
inner automorphisms of $\fk$. This subgroup acts trivially on 
$V(\fk)$ because $\fk$ acts trivially. 
Since all Levi complements are conjugate under 
the group $\Inn(\fk)$ of inner automorphisms 
(cf.\ \cite[Ch.~I]{Bou89}), we obtain with 
$$ \Aut(\fk,\fs) := \{ \phi \in \Aut(\fk) \: \phi(\fs)= \fs\} $$
that 
$$ \Aut(\fk) = \Aut(\fk,\fs) \cdot \Inn(\fk). $$
Since $\Inn(\fk)$ acts trivially on $V(\fk)$, it remains to see that 
the decomposition \eqref{eq:deco} is invariant under $\Aut(\fk,\fs)$. 
Clearly, this group preserves the Levi decomposition of $\fk$, 
hence the subspaces $V(\fs)$ and $V_0$ of $V(\fk)$. 
Moreover, $\Aut(\fs)$ permutes the simple ideals of $\fs$, hence 
preserves the isotypic ideals $\fs_i^{m_i}$ for each $i$. 
This completes the proof. 
\end{proof}

\begin{remark} \label{rem:action-on-V} 
(The action of $\pi_0(\Aut\fs)$ on $V(\fs)$) 
The group $\Aut(\fk)$ acts on the subspace 
$V(\fs)$ on $V(\fk)$ through the natural homomorphism 
$$\Aut(\fk) \to \Aut(\fs),$$ obtained from $\fs \cong  \g/\fr$ 
and the group 
$\Aut (\fs)_{0}$ act trivially on $V(\fs)$ (Remark~\ref{rem:aut-act}). 
The product
$\fS_{m_1} \times \fS_{m_2} \times \ldots \times \fS_{m_r}$ 
of symmetric groups acts naturally on 
$\fs \cong \fs_1^{m_1} \oplus \ldots \oplus \fs_r^{m_r}$
by automorphisms, permuting the simple ideals of $\fs$, 
and since each automorphism of $\fs$ permutes the set of simple ideals, 
we obtain a semidirect decomposition 
$$ \Aut(\fs) \cong \Big(\prod_{i = 1}^r \Aut(\fs_i)^{m_i}\Big) 
\rtimes \prod_{i = 1}^r \fS_{m_i}. $$
This in turn leads to 
$$ \pi_0(\Aut(\fs)) \cong \Big(\prod_{i = 1}^r \pi_0(\Aut(\fs_i))^{m_i}\Big) 
\rtimes \prod_{i = 1}^r \fS_{m_i}. $$

For $V(\s_i)\cong \R$, the invariance of the Cartan--Killing form under 
all automorphisms of $\s_i$ implies that $\Aut(\s_i)$ acts trivially 
on $V(\s_i)$. For $V(\s_i) \cong \C$, the same argument implies that the 
index $2$-subgroup of all complex linear isomorphisms acts trivially 
on $V(\s_i)$, and each antilinear isomorphism $\phi \in \Aut(\s_i)$ acts on 
$V(\s_i) \cong \C$ by complex conjugation. 

If $V(\s_i)$ is one-dimensional, $\fS_{m_i}$ acts by permutations on 
$V(\fs_i^{m_i}) \cong \R^{m_i}$, and 
if $V(\s_i) \cong \C$, then $(\Z/2)^{m_i} \rtimes \fS_{m_i}$ acts by 
permutations on $V(\fs_i^{m_i}) \cong \C^{m_i}$, combined with complex conjugation in the factors.  
\end{remark}

\subsection{The universal period map} 

Let $\kappa_u \: \fk \times \fk \to V(\fk)$ be 
the universal invariant symmetric bilinear 
form. Then $C(\kappa)(x,y,z) := \kappa([x,y],z)$ 
is a $V(\fk)$-valued $3$-cocycle, and the  left invariant closed 
$V(\fk)$-valued $3$-form  
$C(\kappa)^l$ on $K$ specified by $C(\kappa)^l_\1 = C(\kappa)$ defines a period 
homomorphism 
$$ \per_K \:\pi_3(K) \to V(\fk), \quad [\sigma] \mapsto \int_{\sigma} C(\kappa)^l 
= \int_{\bS^3} \sigma^*C(\kappa)^l  $$
(\cite[Lem.~5.7 and Rem.~5.9]{Ne02a}). 
We write $\Pi_K := \im(\per_K)$ for its image. 
To see that this subgroup is fixed by 
$\Aut(\fk)_0 \cong\Aut(K)_0$, we note that for each $\phi \in\Aut(K)$, the
relation 
$$ \per_K \circ \pi_3(\phi) = V(\L(\phi)) \circ \per_K $$
implies that $V(\L(\phi)) \circ \per_K$ only depends on the class 
$[\phi] \in \pi_0(\Aut(K))$. Hence the 
image of $\per_K$ is fixed pointwise by $\Aut(K)_0$. 

\begin{remark} \label{rem:mn.2.3} We recall some results on the homotopy groups of
fi\-nite-di\-men\-sional Lie groups~$K$. 

(a) If $q \: \hat K \to K$ is a covering of Lie groups, then for each 
$j > 1$, the induced homomorphism 
$\pi_j(q) \: \pi_j(\hat K) \to \pi_j(K)$ is an isomorphism. 
This is an easy consequence of the long exact homotopy sequence of the 
principal $\ker q$-bundle $\hat K$ over $K$. 

(b) By E.~Cartan's Theorem, $\pi_2(K) = \1$
(\cite[Th.\ 3.7]{Mim95}). 

(c) Bott's Theorem asserts that for 
a compact connected simple Lie group $K$ we have 
$\pi_3(K) \cong \Z$ (\cite[Th.\ 3.9]{Mim95}). 
A generator of $\pi_3(K)$ can
be obtained from a suitable homomorphism $\eta \: \SU_2(\C) \cong \bS^3\to K$. 
More precisely, let $\alpha$ be a long root in the root system $\Delta_\fk$ 
of $\fk$ and 
$\fk(\alpha) \subeq \fk$ be the corresponding $\su_2(\C)$-subalgebra. 
Then the corresponding homomorphism $\SU_2(\C)\to K$ 
represents a generator of $\pi_3(K)$ (\cite{Bo58}). 
\end{remark}

\begin{remark} \label{rem:pi3} Let 
$K$ be a connected finite-dimensional Lie group,
$C \subeq K$ a maximal compact subgroup, $C_0$ the identity component 
of the center of $C$ and 
$C_1, \ldots, C_m$ be the connected simple normal subgroups of $C$. 
Every compact group is in particular reductive, so that 
the multiplication map 
$$C_0 \times C_1 \times \ldots \times C_m \to C$$
has finite kernel, hence is a covering map.  
As $C_0$ is a torus, its universal covering group is a vector space, 
and therefore 
$\pi_3(C_0) \cong \pi_3(\tilde C_0)$ is trivial. 
Since $K$ is homotopy equivalent to $C$, this leads 
with Remark~\ref{rem:mn.2.3} to 
$$ \pi_3(K) \cong \pi_3(C) \cong \prod_{j = 1}^m \pi_3(C_j) \cong \Z^m. $$
\end{remark}

\begin{proposition} \label{prop:simpcase} Let 
$S$ be a simple connected Lie group with 
Lie algebra~$\fs$. Then 
$$\Pi_S \cong 
  \left\{ 
  \begin{array}{cl} 
\Z & \mbox{for $\fs \not\cong \fsl_2(\R)$} \\ 
\0 & \mbox{for $\fs \cong \fsl_2(\R)$, } 
\end{array} \right. $$
and this group is fixed pointwise by the action of $\Aut(\fs)$ on 
$V(\fs)$. 
\end{proposition}

\begin{proof} Since $\pi_3(S) \cong \pi_3(\tilde S)$ for the universal covering 
Lie group $\tilde S$, we may w.l.o.g.\ assume that $S$ is $1$-connected. 

If $\fs \cong \fsl_2(\R)$, then $S$ is diffeomorphic to $\R^3$,  
so that $\pi_3(S)$ is trivial and therefore $\Pi_S$ is trivial. 
If $\fs \not\cong \fsl_2(\R)$, then the maximal compact subalgebra 
$\fc_\fs$ is not abelian 
(cf.\ \cite[Prop.~VIII.6.2]{Hel78}), so that the maximal compact subgroup 
$C$ of $S$ is non-abelian, hence contains non-trivial simple factors 
$C_1, \ldots, C_m$. In view of Remark~\ref{rem:pi3}, 
$\pi_3(S) \cong \Z^m$ is a non-trivial free group.

For $K := \SU_2(\C)$, pick $x \in \fk$ with 
$\Spec(\ad x) = \{ 0, \pm 2i\}$, where we view $\ad x$ as an endomorphism 
of the complexification $\fk_{\C}\cong \fsl_2(\C)$. The set of all 
such elements is a euclidean 
$2$-sphere in the $3$-dimensional Lie algebra $\su_2(\C)$ which is 
an orbit of the adjoint action. 
Therefore $v_\fk := 4\pi^2 \kappa_u(x,x) \in V(\fk) \cong \R$ 
is well-defined and with Example~\ref{ex:su2} we derive that 
$\Pi_K = \Z  v_\fk.$

Since $\pi_3(S)$ is generated by the homotopy classes of the 
homomorphisms 
$\eta_j \: \SU_2(\C) \to C_j$ specified in Remark~\ref{rem:mn.2.3}(c), 
we conclude that $\Pi_S \subeq V(\fs)$ is the subgroup generated by the 
corresponding elements $v_1,$ $\ldots,$ $v_m$, 
coming from the basis elements $v_{j} =
4\pi^2 \kappa_u(x_j, x_j)\in V(\fc_j),$ where $x_j$ denotes an 
element in a suitable $\su_2$-subalgebra of the simple ideal 
$\fc_j$ of  the maximal compact subalgebra $\fc$ of $\fs$, 
which is normalized in such a way that 
$\Spec(\ad x_j) = \{\pm 2i,0\}$ holds on the $\su_2(\C)$-subalgebra. 
The choice of the elements 
$x_j \in \fc_j$ and the representation theory of 
 $\fsl_2(\C)\cong (\su_2(\C))_{\C}$  imply 
that all eigenvalues of $\ad x_j$ on $\fk_\C$ are 
contained in $i \Z$, so that $\tr ((\ad x_j)^2) \in
-\N_0$. Therefore the values of the Cartan--Killing form of $\fs$ 
on the $x_j$ are integral. 

If $\dim V(\fs)=1$,  then the Cartan--Killing form is universal 
(Remark~\ref{rem:simple}), and  this already implies that 
the elements $v_j$ generate a discrete  non-trivial  subgroup of $V(\fs)$. 
If $\dim V(\fs)=2$, then $\fs$ is complex and $\fc$ is a compact real form 
of $\fs$, hence in particular simple. Therefore  
$\pi_3(S) \cong \pi_3(C) \cong \Z$ (Remark~\ref{rem:mn.2.3}) 
implies that $\Pi_S \cong \Z$. 

To see that $\Aut(\fs)$ fixes $\Pi_S$ pointwise, 
we observe that if $\dim V(\fs) = 1$, then the invariance of the 
Cartan--Killing form under all automorphisms of $\fs$ implies that 
$\Aut(\fs)$ acts trivially on $V(\fs)$. 
If $\dim V(\fs) = 2$, then the subgroup $\Aut_{\C}(\fs)$ of all 
complex linear automorphisms of $\fs$ acts trivially on $V(\fs)$. 
Let $\fc \subeq \fs$ be a compact real form and $\tau \in\Aut(\fs)$ be  
the corresponding antilinear involution. Then 
$\tau$ fixes $\fc$ pointwise, so that the corresponding group 
automorphism fixes $C \subeq S$ pointwise, hence also the canonical image 
$V(\fc) \subeq V(\fs)$, generated by $\Pi_S \cong \Pi_C$. 
Since $\Aut(\fs) \cong \Aut_\C(\fs) \rtimes \{\id,\tau\}$, 
the whole group $\Aut(\fs)$ fixes $\Pi_S$ pointwise. 
\end{proof}

\begin{theorem} \label{thm:disc1} Let $S_i$ be a connected Lie group 
with Lie algebra $\s_i$. Then 
$\Pi_K \cong \prod_{i = 1}^r \Pi_{S_i}^{m_i}$ 
is a discrete subgroup of $V(\fs) \subeq 
V(\fk)$, and if 
$\phi_V \in \GL(V(\fk))$ is induced by an automorphism
$\phi_\fk \in \Aut(\fk)$, then the image of $\Pi_K$ in 
$V(\fk)_{\phi_V}$ is also discrete. 
\end{theorem} 

\begin{proof} We may w.l.o.g.\ assume that $K$ is $1$-connected 
(Remark~\ref{rem:mn.2.3}). 
Then we have a Levi decomposition $K \cong R \rtimes S$, and 
$S \cong S_1^{m_1} \times \ldots \times S_r^{m_r}.$
The functoriality of the period group and the assignment 
$\g \mapsto V(\g)$ now implies that 
$\Pi_K \cong \Pi_S 
\cong \Pi_{S_1}^{m_1} \times \ldots \times\Pi_{S_r}^{m_r},$
where $\Pi_{S_i} \subeq V(\fs_i)$ is a cyclic subgroup, hence discrete 
(Proposition~\ref{prop:simpcase}). 

From Theorem~\ref{thm:decomp} we recall the $\Aut(\fk)$-invariant 
decomposition $V(\fk) = V_0 \oplus \bigoplus_{i=1}^r V_i$ 
with $V_i \cong V(\fs_i)^{m_i}$. We have just seen that the 
period group is adapted to this decomposition with 
$\Pi_{S_i}^{m_i} \subeq V_i$. 
For $i > 0$, $\phi_i := \phi_V\res_{V_i}$ acts on $V_i$ as an element of
$\Aut(\fs_i^{m_i})$, and since $\Aut(\fs_i^{m_i})_0$ acts trivially
(Remark~\ref{rem:aut-act}), 
$\ord(\phi_i) < \infty$ (Remark~\ref{rem:action-on-V}),
so that 
$(V_i)_{\phi_i} \cong V_i^{\phi_i},$
and the projection to the cokernel corresponds to the projection 
to the subspace of fixed vectors of $\phi_i$. 
Since $\phi_i$ preserves $\Pi_{S_i}^{m_i}$ 
(Proposition \ref{prop:simpcase}), 
the image of this group under the projection 
onto the $\phi_i$-fixed space is contained in 
$\frac{1}{\ord(\phi_i)}\cdot \Pi_{S_i}^{m_i}$, hence discrete. 
\end{proof}

\begin{remark} For each $i$, Proposition~\ref{prop:simpcase} implies that the 
the subgroup $\Aut(\fs_i)^{m_i}$ of $\Aut(\fs_i^{m_i})$ fixing 
all simple ideals acts trivially on $\Pi_{S_i^m} \cong \Pi_{S_i}^{m_i}$. 
Therefore only the permutation group $\fS_{m}$ acts 
on this discrete subgroup. 
Thus any automorphism of $\fk$ acts on $\Pi_K 
\cong \prod_{i = 1}^r \Pi_{S_i}^{m_i}$ as an element 
of the product $\fS_{m_1} \times \ldots \times \fS_{m_r}$. 
\end{remark}

\section{Appendix: Central extensions of Lie groups} 
\label{app:C}

In this appendix we recall some facts on the integration of Lie algebra 
$2$-cocycles from \cite{Ne02a}. 
They provide a general set of tools to integrate 
central extensions of Lie algebras to extensions of connected Lie groups.

Let $G$ be a connected Lie group and $V$ be a Mackey complete space. 
Further, let $\omega \in Z^k(\g,V)$ be a $k$-cocycle and 
$\omega^l \in \Omega^k(G,V)$ be the corresponding left equivariant 
$V$-valued $k$-form with $\omega^l_\1 = \omega$. 
Then each continuous map $\bS^k \to G$ is homotopic to a smooth map 
(cf.\ \cite{approx} or \cite{Ne02a}), and 
$$ \per_\omega \: \pi_k(G) \to V, \quad 
[\sigma] \mapsto \int_\sigma \omega^l 
= \int_{\bS^k} \sigma^*\omega^l 
= \int_{\bS^k} \omega(\delta\sigma,\ldots, \delta\sigma), $$
for $\sigma \in C^\infty(\bS^k,G)$, defines the {\it period homomorphism} 
whose values lie in the $G$-fixed part of $V$ 
(\cite[Lem.~5.7 and Rem.~5.9]{Ne02a}). 

\begin{theorem} {\rm(\cite[Prop.~7.6 and Thm.~7.9]{Ne02a})}\label{thm:int-general} 
Let $G$ be a connected Lie group with Lie algebra
$\g$. A central Lie algebra extension 
$\hat \g = V \oplus_\omega \g$ defined by $\omega \in Z^2(\g,V)$ 
integrates to a Lie group extension of some covering group 
of $G$ if and only if the period group 
$\Pi_\omega := \per_\omega(\pi_2(G)) \subeq V$ is discrete. 
It integrates to an extension of $G$ if and only
if the adjoint action of $G$ on $\g$ 
lifts to an action of $G$ on $V \oplus_\omega \g$.
\end{theorem}

\begin{remark} \label{rem:8.4} (a) To 
calculate period homomorphisms, it is often 
convenient to use related cocycles on different groups. So, let 
us consider a morphism $\phi \: G_1 \to G_2$ of Lie groups and 
$\omega_2 \in Z^2(\g_2, V)$, $V$ a trivial $G_i$-module. 
Then a straight forward argument 
shows that 
\begin{equation}
  \label{eq: per-trans}
\per_{\omega_2} \circ \pi_2(\phi) = \per_{\L(\phi)^*\omega_2} 
\: \pi_2(G_1) \to V. 
\end{equation}

\nin (b) From (a) we obtain in particular 
for $\omega \in Z^2(\g,V)$ and $\phi \in \Aut(G)$ the relation 
\begin{equation}
  \label{eq: per-auto}
\per_{\omega} \circ \pi_2(\phi) = \per_{\L(\phi)^*\omega} 
\: \pi_2(G) \to V. 
\end{equation}
If, in addition, $\phi$ is homotopic to the identity in 
the sense that $\phi = \gamma(1)$ for a curve 
$\gamma \: [0,1] \to G$ with $\gamma(0)= \id_G$ for which the map 
$$ \tilde\gamma \: [0,1] \times G \to G, \quad (t,g) \mapsto 
\gamma(t)(g) $$ 
is smooth, then $\pi_2(\phi)= \id$ implies that the periods  
of the $2$-cocycle 
$$\L(\phi)^*\omega - \omega$$
are trivial. In this case we further have a derivation 
$D = \phi'(0) \in \der(\g)$, and by applying the preceding relation 
to all automorphisms and taking derivatives in $\1$, it follows 
that the periods of the cocycle 
$$ \omega_D(x,y) = \omega(Dx,y) + \omega(x,Dy) 
= \frac{d}{dt}\Big|_{t=0} \L(\phi_t)^*\omega $$
vanish. 
\end{remark}

The following theorem can be found in \cite[Thm.~V.9]{MN03}: 

\begin{theorem} \label{thm:lifting-theorem} {\rm(Lifting Theorem)}  
Let $q \: \hat G \to G$ be a central Lie group extension of the $1$-connected 
Lie group $G$ by the Lie group $Z \cong \z/\Gamma_Z$. 
Let $\sigma_G \: H \times G \to G$, resp., $\sigma_Z \: H \times Z \to Z$ 
be smooth automorphic actions of the 
Lie group $H$ on $G$, resp., $Z$ and 
$\sigma_{\hat\g}$ be a smooth action of $H$ on $\hat\g$ compatible with the 
actions on $\z$ and $\g$. Then there is a unique smooth action 
$\sigma_{\hat G} \: H \times \hat G \to \hat G$ 
by automorphisms compatible with the actions on $Z$ and $G$, for which 
the corresponding action on the Lie algebra $\hat\g$ 
is $\sigma_{\hat\g}$. 
\end{theorem}

\section{Appendix: Some facts on curvature and parallel transport} 
\label{app:D}

Let $M$ be a finite-dimensional manifold, $H$ a regular Lie group, 
$P(M,H,q)$ a principal $H$-bundle over $M$ and $\theta \in \Omega^1(P,\fh)$ 
a principal connection $1$-form. 

For each piecewise smooth curve 
$\alpha \: [a,b] \to M$, we then have an $H$-equivariant 
{\it parallel transport map}
$\Pt(\alpha) \: P_{\alpha(a)} \to P_{\alpha(b)}$ defined by 
$\Pt(\alpha).p := \hat\alpha(1)$, 
where $\hat\alpha \: [a,b] \to P$ is the horizontal 
lift of $\alpha$ starting in $p$. For a closed curve 
parallel transport and holonomy are connected by 
$$ \Pt(\alpha).p_0 = p_0.{\cal H}(\alpha). $$

Let $v,w \in T_{m_0}(M)$ and 
consider an open connected neighborhood $U$ of $m_0$ in $M$ such that 
$P\res_U$ is trivial and there exist smooth vector fields 
$X,Y \in {\cal V}(U)$ with $X(m_0) = v$ and $Y(m_0) = w$ and 
$T > 0$ such that for $0 \leq t_i \leq T$ the points 
$$ \Fl^Y_{-t_4} \circ \Fl^X_{-t_3} 
\circ \Fl^Y_{t_2} \circ \Fl^X_{t_1}(m_0) $$
are defined and contained in $U$. 

For $0 \leq t \leq T$, 
$$\gamma(t) := \Fl^Y_{-t} \circ \Fl^X_{-t} \circ \Fl^Y_{t} 
\circ \Fl^X_{t}(m_0). $$
defines a smooth curve in $M$ with $\gamma(0) = m_0$, $\gamma'(0) = 0$, 
and 
$$ \gamma''(0) = 2 [X,Y](m_0). $$
(\cite[Thm.~1.4.4]{BC64}).  

We write $\alpha_t \: [0,5t] \to M$ for the curve obtained 
by concatenating integral curves of $X$, $Y$, $-X$ and $-Y$ 
defined on $[0,t]$ with the reversed curve $\gamma$, so that we obtain 
a loop in $m_0$ which is piecewise smooth. 
Note that any piecewise smooth loop can be reparametrized as a smooth loop, 
so that $\beta(t) := {\cal H}(\alpha_t)$ also is the holonomy of a smooth 
loop, and it is clear that it is a smooth curve in $H$. We claim that 
$$ \beta'(0) = 0 
\quad \mbox{ and } \quad 
\beta''(0) = 2 R(\theta)_{p_0}(\tilde v,\tilde w), $$
where $\tilde v \in T_{p_0}(P)$ denotes the unique horizontal lift 
of $v \in T_{m_0}(M)$ 
(cf.\ \cite[Thm.~6.1.3]{BC64}). 

Since the bundle $P_U$ is trivial, 
we may w.l.o.g.\ assume that $P_U = U \times H$. Then the 
connection $1$-form $\theta$ has the form 
$$ \theta = p_H^*\kappa_H + \Ad(p_H)^{-1}(p_U^*A), $$
where $A \in \Omega^1(U,\fh)$, and $p_H \: U \times H \to H$ and 
$p_U \: U \times H \to U$ are the projection maps. After adjusting the 
trivialization if necessary, we may w.l.o.g.\ assume that 
$A(m_0) = 0$, i.e., the subspace $T_{m_0}(M)$ is horizontal 
in $T_{(m_0,\1)}(P)$. 

Let $\tilde X, \tilde Y \in {\cal V}(P)^H$ denote the unique 
horizontal lifts of the vector fields $X, Y$ and 
$$ \hat\gamma(t) := \Fl^{\tilde Y}_{-t} \circ \Fl^{\tilde X}_{-t} 
\circ \Fl^{\tilde Y}_{t} \circ \Fl^{\tilde X}_{t}(m_0), $$
which coincides with $\hat\alpha_t(4t)$ for the horizontal lift 
$\hat\alpha$ of $\alpha$, starting in $p_0$. 
In the product coordinates of $P_U = U \times H$, we now find with 
$p_0 = (m_0,\1)$: 
$$ \hat\gamma(t) = (\gamma(t), \zeta(t)), $$
where 
$\zeta(0) = \zeta'(0) = 0$ and 
$$ \zeta''(0) 
= 2 \theta([\tilde X, \tilde Y])(p_0)
= 2 \dd\theta(\tilde Y, \tilde X)(p_0) 
= 2 R(\theta)(\tilde w, \tilde v). $$
Let $(\gamma(s),\rho(s))$, $0 \leq s \leq t$, denote the horizontal lift 
of the curve $\gamma$ starting in $p_0$. Then 
$(\gamma,\rho).\beta(t) = (\gamma,\rho\cdot \beta(t))$ 
is the horizontal lift starting in 
$p_0.\beta(t)$, and we thus obtain 
$$\rho(t)\cdot \beta(t) = \zeta(t). $$

Since $(\gamma,\rho)$ is horizontal, we have 
$\delta(\rho)_t = - A_{\gamma(t)}(\gamma'(t)),$
which leads to 
$\rho(0) = \1, \rho'(0) = 0$
and further to  
$$ \rho''(0) = - A_{m_0}(\gamma''(t)) = 0. $$
Hence 
$$ \beta''(0) = (\rho \cdot \beta)''(0) = 
\zeta''(0) = 2 R(\theta)(\tilde w, \tilde v). $$

We have thus constructed a family $(\alpha_t)_{0 \leq t \leq T}$ of 
(piecewise) smooth loops in $m_0$ for which the holonomy defines a smooth 
curve $\beta(t) = {\cal H}(\alpha_t)$ in $H$ with 
$\beta'(0) = 0$ and $\beta''(0) = 2 R(\theta)(\tilde w,\tilde v)$.

\end{document}